\documentclass[10pt]{amsart} 
\usepackage{amsmath, amsfonts, amsthm, amssymb, amscd}
\usepackage[mathcal]{euscript}  
\usepackage{dsfont, mathrsfs, accents, verbatim} 
\theoremstyle{plain}
\newtheorem{theorem}{Theorem}[section]
 
\newtheorem{lemma}[theorem]{Lemma}
\newtheorem{corollary}[theorem]{Corollary} 
\newtheorem{definition}[theorem]{Definition} 
   
\numberwithin{equation}{section} 
\newcommand \RR   {\mathbb R}

\newcommand \be  {\begin{equation}}
\newcommand \ee  {\end{equation}}
\newcommand \bs  {\begin{split}}
\newcommand \es  {\end{split}} 
\newcommand \del 	\partial
\newcommand \eps \varepsilon
\newcommand \lam \lambda
 
\newcommand \Rcal {\mathcal R} 
\newcommand \Ucal {\mathcal U} 
\newcommand \Kcal {\mathcal K} 
\newcommand \Fcal {\mathcal F} 
\newcommand \Tcal {\mathcal T} 
   
\newcommand \Rbf {\mathbf R} 
\newcommand \Sbf {\mathbf S} 
\newcommand \ubar {\overline u} 
\newcommand \vbar {\overline v} 
\newcommand \abar {\overline a} 
\newcommand \Rbar {\overline R} 
\newcommand \bbar {\overline b}

\newcommand \sgn {{\rm sgn }\,} 
\newcommand \one {{\mathds 1}} 
\newcommand \auth \textsc
\newcommand \la {\langle}
\newcommand \ra {\rangle}
 
\newcommand \sigmas {\sigma^\sharp}
\newcommand \supp {\text{supp }} 

\begin{document}

\bibliographystyle{plain}

\title[Entropy solutions of the Euler equations]
{Entropy solutions of the Euler equations
\\ for isothermal relativistic fluids}  
\author[P.G. LeFloch and M. Yamazaki]
%
{\scshape Philippe G. LeFloch$^1$ and Mitsuru Yamazaki$^2$}
\thanks
{$^1$ Laboratoire Jacques-Louis Lions 
\& Centre National de la Recherche Scientifique, Universit\'e de Paris 6, 
4 place Jussieu, 75252 Paris, France. 
E-mail : {\sl lefloch@ann.jussieu.fr}}
\thanks{$^2$ Graduate School of Pure and Applied 
Sciences, University of Tsukuba, 305-8571
Ibaraki, Japan. 
E-mail : {\sl yamazaki@math.tsukuba.ac.jp}
\newline
2000\textit
}

\maketitle

\subsection*{Abstract}{\small 
We investigate the initial-value problem for the relativistic Euler equations governing isothermal perfect fluid flows,
and generalize an approach introduced by LeFloch and Shelukhin in the non-relativistic setting.  
We establish the existence of globally defined, bounded measurable, entropy solutions with arbitrary large
amplitude. An earlier result by Smoller and Temple for the same system covered solutions with 
bounded variation that avoid the vacuum state. The new framework proposed here 
provides entropy solutions in a larger function space and allows for the mass density  
to vanish and the velocity field to approach the speed of light. 
The relativistic Euler equations become strongly degenerate in both regimes, as the conservative 
or the flux variables vanish or blow-up. Our proof is based on the method of compensated compactness 
for nonlinear systems of conservation laws (Tartar, DiPerna) and 
takes advantage of a scaling invariance property 
of the isothermal fluid equations. We also rely on properties of the fundamental kernel 
that generates the mathematical entropy and entropy flux pairs. This kernel exhibits certain singularities
on the boundary of its support and we are led to analyze certain nonconservative products 
(after Dal~Maso, LeFloch, and Murat) consisting of functions of bounded variation by measures. 
} 
 

\newpage 


\section{Introduction and main result}
\label{IN-0} 

\subsection{Purpose of this paper}

The relativistic Euler equations describe the dynamics of a compressible fluid in the context 
of special relativity, i.e.~a fluid evolving on the flat Minkowski spacetime. 
This system can be regarded as an approximation to the Euler-Einstein equations, which is 
valid in a small region of the spacetime and away from large matter concentrations. 
Under the assumption of plane-symmetry, the fluid unknowns consist of the co-moving mass density ${\rho \geq 0}$ 
and the velocity field ${v \in (-c,c)}$, where $c$ denotes the speed of light. 
Our main purpose in the present paper is the construction of weak solutions containing 
shock waves and having arbitrary large amplitude. 

The mathematical properties of this model were investigated in works 
by Taub (1957), Thompson (1986), Lichnerowicz (1993), and others. 
These equations are also important in computational physics 
and we refer to Marti and M\"uller (2003) for an extensive review of the Riemann problem and 
the numerical methods in hydrodynamics in the context of special relativity. 
The mathematical analysis of the relativistic Euler equations for the isothermal fluids,
considered in the present paper, has received much less attention in the literature; 
our emphasis is on fluids governed by the linear pressure law 
$$
p(\rho) = k^2 \, \rho,
$$
where $k>0$ is the (local, constant) sound speed and by the principle of special relativity 
must be less than the speed of light denoted here by $c$.

The Euler equations form a nonlinear system of partial differential equations of hyperbolic type. 
It is well-known that solutions even they are smooth initially will eventually become discontinuous
and must be understood in the weak sense of distributions. Furthermore, for the sake of solutions
it is necessary to constrain these weak solutions by certain entropy inequalities. 

The initial-valued problem for relativistic isothermal fluids 
was first studied by Smoller and Temple (1993). They  
established the existence of entropy solutions, under the 
assumption that the initial mass density is bounded, is bounded away from zero, and has bounded variation. 
Their result is based on the Glimm scheme and extends an earlier approach 
for the non-relativistic version (Nishida, 1968). 

In the present paper, we propose an alternative approach based on the 
method of compensated compactness 
for nonlinear conservation laws (Tartar, 1983) 
and, for the relativistic Euler equations,  
we provide a mathematical framework encompassing a large class of weak solutions. These solutions may take vacuum values ${\rho = 0}$ and may have high-velocity
approaching the speed of light. Solutions of this nature arise naturally in applications;   
for instance, a star is described by a compactly supported mass density function.  
As can be checked easily, at points where $\rho$ vanishes the Euler equations are highly degenerate 
and the conservative variables vanish identically while the velocity field is ill-defined. 
Another related singularity of the equations is obtained in the limit when the fluid velocity approaches 
the light speed and the wave speeds approach each other; hence, the model fails to be uniformly strictly hyperbolic. 
These features of the relativistic Euler equations for isothermal fluids 
lead to particularly challenging mathematical questions, 
concerning the existence and the behavior of entropy solutions. 


\subsection{Relativistic fluid equations}

We consider the following system of two conservation laws
\be
\nonumber 
\aligned
&  \del_t\left(\frac{\rho \, c^2+ p \, v^2/c^2}{ c^2-v^2} \right)
     + \del_x \left( \frac{ (p+\rho c^2) \, v }{ c^2 - v^2 } \right) = 0,\\
&   \del_t \left( \frac{(p + \rho \, c^2)  \, v }{ c^2 - v^2} \right) 
  + \del_x \left( \frac{(p + \rho \, v^2) \, c^2 }{ c^2 - v^2 } \right) = 0,   
\endaligned
\ee 
which, by setting $\eps = 1/c$ and using the condition $p(\rho) = k^2 \, \rho$, reads   
\be
\label{IN.Euler2}
\begin{aligned}  
&   \del_t\left( \frac{1 + \eps^4 \, k^2 \, v^2 }{ 1- \eps^2 \, v^2} \, \rho\right)
  + \del_x \left(  \frac{1 + \eps^2 \, k^2 }{ 1 - \eps^2 \, v^2} \, \rho \, v  \right) = 0,
\\
&   \del_t \left(  \frac{ 1 + \eps^2 \, k^2 }{ 1 - \eps^2 \, v^2} \, \rho \, v  \right)
  + \del_x \left( \frac{k^2 +  v^2}{ 1 - \eps^2 \, v^2} \, \rho \right) = 0. 
\end{aligned}
\ee
Taking the formal limit $\eps=0$, we obtain the non-relativistic version of these equations: 
\be
\label{IN.Euler-non}
\begin{aligned}  
& \del_t \rho  + \del_x \left( \rho \, v  \right) = 0,
\\
& \del_t \left( \rho \, v  \right) + \del_x \left( (k^2 +  v^2) \, \rho  \right) = 0. 
\end{aligned} 
\ee

Our starting point is the recent work by Huang and Wang (2003)  
and LeFloch and Shelukhin (2005). 
The existence of entropy solutions to \eqref{IN.Euler-non} is known when the mass density is bounded
and the velocity is unbounded. It was observed that the natural function space associated with  
the non-relativistic Euler equations allows for the velocity field to be unbounded. 
Our aim will be here to generalize to relativistic fluids the method and results by 
LeFloch and Shelukhin (2005), 
based on such ``natural'' estimates.

The equations \eqref{IN.Euler2} form a nonlinear hyperbolic system of partial 
differential equations of the form  
\be
\label{IN.cons} 
\begin{aligned}
& \del_t G + \del_x H = 0, \qquad \del_t H + \del_x F = 0, 
\end{aligned}
\ee
where 
$$
\begin{aligned}
& G(\rho,v)  =  \frac{1 + \eps^4 k^2 v^2}{1 - \eps^2 v^2} \rho,
\quad H(\rho,v) =  \frac{1 + \eps^2 k^2}{1 - \eps^2 v^2} \rho v,
\\
& F(\rho,v) =  \frac{k^2 + v^2}{1 - \eps^2 v^2} \rho. 
\end{aligned}
$$  
Since shock waves are known to arise even from smooth initial data (Pan and Smoller, 2006), 
we need a concept of solutions that include discontinuous functions. 

We introduce the following notion of entropy solution.

First of all, we will say that a Lipschitz continuous map $(\Ucal,\Fcal)$ is an {\sl entropy pair} 
if every {\sl smooth} solution of \eqref{IN.Euler2} satisfies the additional conservation law
$$
\del_t \Ucal(\rho,v) + \del_x \Fcal(\rho,v) = 0. 
$$
However, weak solutions are required to satisfy the above conservation laws for convex functions $\Ucal$,  
but as inequalities only. 
Furthermore, the class of entropy functions is further restricted as we only consider 
{\sl weak entropies} that is functions $\Ucal$ vanishing on the vacuum line $\rho =0$. 
For instance, both pairs $(G,H)$ and $(H,F)$ in \eqref{IN.cons} are (trivial) weak entropy pairs.  
We also set  
$$
\eps' := \frac{2 \eps }{ 1+\eps^2} \in (0,1). 
$$

\begin{definition}[Notion of entropy solutions] 
\label{AS-tame}
A  {\rm tame region} is a set of the form 
$$
\Tcal_\eps(M) := \big\{ \rho,v \, / \,  0 \leq \rho \leq M, \quad 1 - \eps |v|  \geq (\rho/M)^{\eps'} \big\}
$$ 
for some constant $M>0$. A pair of measurable and bounded functions $\rho_0,v_0 : \RR \to \RR$ is a {\rm tame initial data} 
if its range is included in a tame region. 

Given a tame initial data $\rho_0, v_0$, a pair of measurable and bounded
 functions $\rho,v :\RR_+ \times \RR \to \RR$ is called a {\rm tame entropy solution} 
to the isothermal relativistic Euler equations \eqref{IN.Euler2} associated with the initial data $\rho_0, v_0$ 
if the range of $\rho,v$ is included in a tame region and, moreover, 
$$
\iint_{\RR_+ \times \RR} \Big( \Ucal(\rho,v) \, \del_t \theta + \Fcal(\rho,v) \, \del_x \theta \Big) \, dxdt 
+ \int_\RR \Ucal(\rho_0, v_0) \, \theta(0, \cdot) \, dx \geq 0 
$$
for every convex, smooth, weak entropy pair $(\Ucal,\Fcal)$ of the isothermal relativistic Euler equations
and every non-negative test-function $\theta$ supported in $[0,\infty) \times \RR$. 
\end{definition}

Observe that the inequality in the definition of a tame region allows the velocity to approach the light speed 
(normalized here to be $1/\eps$) when the mass density approaches zero. The concept of
 a tame region is quite natural, as it is equivalent to uniform bounds on the Riemann invariants 
(defined later in Section~\ref{BP-0}). 
Clearly, the entropy pairs under consideration need not be globally Lipschitz continuous, 
but only so within any tame region, 
as those are the only regions of interest. 

One of our main results in this paper is the following existence theorem.

\begin{theorem}[Existence theory] 
\label{AS-maintheorem} 
Given any parameter value $\eps \in (0,1)$ and a tame initial data $\rho_0, v_0$, the initial-value problem 
for the relativistic Euler equations for isothermal fluids \eqref{IN.Euler2} 
admits a tame entropy solution $\rho,v :\RR_+ \times \RR \to \RR$ associated with $\rho_0, v_0$. 
\end{theorem}


\subsection{Main ideas for the proof} 

The approach proposed in the present paper will rely on the following observation, which was already 
pointed out in the non-relativistic setting.

\begin{lemma}[Linearity property] 
\label{invariance}
If $(\rho,v)$ is a (weak, entropy) solution of the relativistic Euler equations for isothermal fluid \eqref{IN.Euler2}, 
then for every positive constant $\lambda$, the function $(\lambda \rho, v)$ is also a 
(weak, entropy) solution of the same equations.
\end{lemma}

Our general strategy of proof follows, on one hand, DiPerna (1983), who obtained bounded solutions
$\rho \geq 0$ and $v \in \RR$ for non-relativistic polytropic fluids 
satisfying $p(\rho) = k^2 \rho^\gamma$, with $\gamma>1$ and, on the other hand,  
LeFloch and Shelukhin (2005), who extended DiPerna's analysis to include isothermal fluids satisfying 
$p(\rho) = k^2 \, \rho$ and observed that the velocity field $v$ need not be bounded.

The main difficulty for our analysis in this paper lies in the lack of uniform 
strict hyperbolicity of the Euler equations when the fluid velocity approaches the light speed. 
To deal with this problem we will proceed along the following lines: 
\begin{itemize}
\item {\sl Mathematical entropy pairs.} Our first task will be 
to construct entropy pairs which amounts to solve a linear hyperbolic equation in the variable $\rho,v$. 
We will introduce the Riemann function $\Rcal$ and the entropy kernel $\chi$ associated with this equation, 
so that the entropy pairs of interest can be expressed by an explicit formula in term of the kernel $\chi$. 
Contrary to the case of polytropic fluids (DiPerna, 1983) 
and in agreement with the case of (non-relativistic) isothermal fluids, 
an initial data for the entropy kernel must be imposed away from the vacuum, say on the line $\rho=1$. 
\item {\sl Singularities of the entropy kernel.} 
We will show that the function $\chi$ is discontinuous along the boundary $\del \Kcal$ 
of its support, so that its first-order derivatives exhibit Dirac masses. 
In the case $\eps =0$ the kernel was given by an explicit formula. 
In contrast, when $\eps >0$ we need to derive uniform estimates on $\chi$ and determine 
explicitly the traces of its first-order derivatives along $\del\Kcal$.  
\item {\sl A~priori bounds.} We will  next derive a~priori bounds on approximate solutions $\rho^h, v^h$ 
generated by the Lax-Friedrichs scheme. Our bounds show that the mass density remains uniformly bounded 
and the velocity field satisfies the tame condition.  
\item{\sl Reduction of the Young measure.} In the final part of the proof we identify the structure 
of a Young measure $\nu=\nu_{t,x}$ associated with the sequence $\rho^h, v^h$. We
 analyze certain nonconservative products 
(Dal Maso, LeFloch, and Murat, 1995) consisting of functions of bounded variation by measures.
The term of interest contains a key coefficient (denoted below by $\Xi(\rho)$) which does not vanish, 
provided we take advantage of the scaling invariance property in Lemma~\ref{invariance}. 
\end{itemize}


\section{Basic properties of the model}
\label{BP-0}

\subsection{Wave speeds and Riemann invariants}

Scaling properties of the equations \eqref{IN.Euler2} will play an important role. Observe that the 
transformation  
$$
v' = v/k, \quad t'= k \, t, \quad \eps' = k \, \eps
$$
allows one to reduce the system \eqref{IN.Euler2} to the same system with $k=1$. In view of the physical constraint 
$0 < k < c$ between the sound speed and the light speed this amounts to impose 
$c>1$. The limiting case $c \to 1$ corresponds to the special case where the sound speed and the 
light speed coincide. From now on, we suppose that $k=1$ so that the Euler equations read 
\be\label{BP.Euler3} 
\begin{aligned}  
&   \del_t\left( \frac{1 + \eps^4 \,  v^2 }{ 1- \eps^2 \, v^2} \, \rho\right)
        + \del_x \left(  \frac{1 + \eps^2  }{ 1 - \eps^2 \, v^2} \, \rho \, v  \right) = 0,
\\
&   \del_t \left(  \frac{ 1 + \eps^2  }{ 1 - \eps^2 \, v^2} \, \rho \, v  \right)
       + \del_x \left( \frac{1 +  v^2}{ 1 - \eps^2 \, v^2} \, \rho  \right) = 0. 
\end{aligned}
\ee 
The velocity is restricted to lie in the interval $(-1/\eps,1/\eps)$; note that the conservative and flux variables 
blow-up when $v \to \pm 1/\eps$. The range of physical interest for $\eps$ is
$$
0 < \eps < 1, 
$$
the limiting case $\eps =0$ and $\eps=1$ corresponding to the non-relativistic model (speed of light is infinite) 
and the scalar field model (the sound speed and the light speed coincide), respectively. 

Indeed, the system \eqref{BP.Euler3} in the limit $\eps \to 1$ converges to 
\be
\label{scalarfield} 
\begin{aligned}  
&   \del_t\left( \frac{1 + v^2 }{ 1- v^2} \, \rho\right)
        + \del_x \left(  \frac{2}{ 1 - v^2} \, \rho \, v  \right) = 0,
\\
&   \del_t \left(  \frac{2}{ 1 - v^2} \, \rho \, v  \right)
       + \del_x \left( \frac{1 +  v^2}{ 1 - v^2} \, \rho  \right) = 0,
\end{aligned}
\ee 
which is simply equivalent to the linear wave equation. 
This is clear by introducing the unknowns 
$a := \frac{1 + v^2 }{ 1- v^2} \, \rho$ and $b :=  \frac{2}{ 1 - v^2} \, \rho \, v$, 
so that
\be
\nonumber 
\begin{aligned}  
&   \del_t a + \del_x b = 0,
\\
&   \del_t b + \del_x a = 0. 
\end{aligned}
\ee 

\

The conservation laws \eqref{BP.Euler3} form a nonlinear hyperbolic system whose Jacobian matrix 
admits the two eigenvalues 
$$ 
\lam_1 := \frac{ v - 1 }{ 1 - \eps^2 \, v}, \qquad \lam_2 := \frac{ v + 1 }{ 1 + \eps^2 \, v}. 
$$
Clearly, the characteristic speeds are smooth functions in the closed interval $v \in [-1/\eps,1/\eps]$.  
The corresponding eigenvectors are 
$$
\aligned 
& r_1 := \Big( \frac{-1 }{ 1 - \eps^2 \, v^2}, \, \frac{1 }{ 1 + \eps^2}\, \frac{1 }{ \rho}  \Big), 
\\
& r_2 = \Big( \frac{1 }{ 1 - \eps^2 \, v^2}, \, \frac{1 }{ 1 + \eps^2}\, \frac{1 }{ \rho}  \Big).
\endaligned 
$$

The Riemann invariants $w,z$, by definition, satisfy 
$\nabla w \cdot r_1 = 0$, $\nabla z \cdot r_2 = 0$, 
and are uniquely defined up to the composition by a one-to-one map: 
\be
\label{Rinvariants}
\aligned
& w := u + R,       & R= \frac{w - z }{ 2}, \\
& z := u - R,        &  u = \frac{w + z }{ 2}, 
\endaligned
\ee
where $R,u$ are functions of $\rho,v$:  
$$ 
\begin{aligned}  
&   R = \Rbar(\rho) : = \frac{1 }{ 1 + \eps^2} \, \ln \rho, \,  &  \rho  =&e^{(1 + \eps^2) R},
\\
&   u = \ubar(v) := \frac{1 }{ 2 \, \eps} \, \ln\Big( \frac{1 + \eps \, v }{ 1 - \eps \, v} \Big), \,   & v  =& \frac1\eps\frac{e^{2\eps u} - 1}{e^{2\eps u} + 1}. 
\end{aligned} 
$$
The Riemann invariants provide a change of variables $(\rho,v) \mapsto (w,z)$, 
which will be often used. 

Clearly, the mapping $v \mapsto \ubar(v)$ is one-to-one from 
the bounded interval $(-1/\eps,1/\eps)$ onto the real line $\RR$. 
The mapping $\rho \mapsto R(\rho)$ is one-to-one from $(0, \infty)$ onto 
$\RR$. It is not difficult to check that, 
in terms of the variables $w,z$, the system \eqref{BP.Euler3} takes the diagonal form
$$
\begin{aligned}  
& \del_t w + \lam_2 \, \del_x w = 0,
\qquad  \del_t z + \lam_1 \, \del_x z = 0.
\end{aligned} 
$$ 
Observe that 
\be
\label{speeds}
\aligned 
& \lambda_1(w,z) =  -\frac1{\eps} \frac{1 + \eps - (1-\eps) e^{\eps (w+z)}}{1+ \eps + (1-\eps) e^{\eps (w+z)}},  
\\
& \lambda_2(w,z) = -\frac1{\eps} \frac{1 - \eps - (1+\eps) e^{\eps (w+z)}}{1- \eps + (1+\eps) e^{\eps (w+z)}} = - \lambda_1 (-w,-z). 
\endaligned 
\ee 

Sometimes, we will also use of the ``modified'' Riemann invariants defined as
\be\label{modifiedRinvariants}
\aligned 
& W := e^w =  \rho^{1/(1+\eps^2)}  \Big({1 + \eps \, v \over 1 - \eps \, v}\Big)^{1/(2\eps)}, 
\\ 
& Z := e^{-z} =  \rho^{1/(1+\eps^2)}  \Big({1 + \eps \, v \over 1 - \eps \, v}\Big)^{-1/(2\eps)}. 
\endaligned 
\ee
Note that $\rho = (WZ)^{(1 + \eps^2)/2}$. 

Expressing now the physical variables $\rho, v$ as functions of the Riemann invariants, 
$$
\rho = \Rbar^{-1}\left( (w-z)/2 \right) = \exp\left((1 + \eps^2)\frac{w - z}2\right),
$$
and
$$
\begin{aligned}  
  v = \vbar (w + z) : 
 & =   \frac{1 }{ \eps}  \,  \frac{ e^{\eps (w+z) } - 1 }{ e^{ \eps (w+z) } + 1 }\\
 &   = \frac{1 }{ \eps} \, \Big( 1- \frac{2 }{ e^{ \eps (w + z)} + 1 } \Big) 
  = \frac1\eps  \tanh \Big(\eps \frac{w + z}2 \Big), 
\end{aligned} 
$$
we obtain the $w$- and $z$-derivatives of $(\rho,v)$: 
$$
\aligned 
& v_w = v_z = \frac1{2 u_v} = \frac{1 }{ 2} \, \Big(1 - \eps^2 \, v^2 \Big),
\\
& \rho_w = -\rho_z = \frac{1 }{ 2 \, R_\rho} = \frac{ 1 + \eps^2  }{ 2 }\rho.
\endaligned 
$$
Moreover, the derivatives of the Riemann invariants considered as functions of $(\rho, v)$ are 
$$
\aligned 
& w_\rho = -z_\rho = \frac{1 }{ 1 + \eps^2}\, \frac{1 }{ \rho}, 
\\
& w_v = z_v = \frac{1 }{ 1-\eps^2 \, v^2}. 
\endaligned 
$$

Finally, for latter use we express the Euler equations in the nonconservative variables $(\rho,v)$. 
By setting 
\be
\nonumber 
\aligned
{D(G,H) \over D(\rho,v)} = 
\left(
\begin{array}{cc}
 \displaystyle\frac{1+\eps^4 v^2}{1 - \eps^2 v^2} & \displaystyle\frac{2 \eps^2 (1 + \eps^2)\rho v}{(1 - \eps^2 v^2)^2}\\
 \displaystyle\frac{(1 + \eps^2) v}{1 - \eps^2 v^2} & \displaystyle\frac{(1 + \eps^2 ) \rho (1 + \eps^2 v^2)}{(1 - \eps^2 v^2)^2}
\end{array}
\right) 
\endaligned
\ee
and 
\be
\nonumber 
\aligned
{D(H,F) \over D(\rho,v)} = \left(
\begin{array}{cc}
 \displaystyle\frac{(1 + \eps^2 )v}{1 - \eps^2 v^2} & \displaystyle\frac{(1 + \eps^2 ) \rho (1 + \eps^2 v^2)}{(1 - \eps^2
v^2)^2}\\  
\displaystyle\frac{1 + v^2}{1 - \eps^2 v^2} & \displaystyle\frac{2 (1 + \eps^2) \rho v}{(1 - \eps^2 v^2)^2}
\end{array}
\right),
\endaligned
\ee
we can rewrite \eqref{BP.Euler3} as 
$$ 
\del_t \widetilde u + \del_x \widetilde G(\widetilde u) = 0, 
\qquad \widetilde u := \left(\begin{array}{c}\rho\\v\end{array}\right), 
$$
with 
$$ 
\aligned 
{D\widetilde G \over D\widetilde u}  
& = \left({D(G,H) \over D(\rho,v)}\right)^{-1} {D(H,F) \over D(\rho,v)}
\\
& = \left(
\begin{array}{cc}
\displaystyle\frac{(1 - \eps^2 )v}{1 - \eps^4 v^2} & \displaystyle\frac{(1 + \eps^2 ) \rho}{1 - \eps^4 v^2}\\
\displaystyle\frac{(1 - \eps^2 v^2)^2}{(1 + \eps^2) \rho (1 - \eps^4 v^2)} & \displaystyle\frac{(1 - \eps^2) v}{1 - \eps^4 v^2}
\end{array}
\right).
\endaligned 
$$


\subsection{Strict hyperbolicity fails for high-velocity fields} 

From the relation  
$$
\lam_2 - \lam_1 = 2 \, \frac{1 - \eps^2 v^2 }{ 1 - \eps^4 v^2} > 0 
$$
we deduce that: 

\begin{lemma}[Hyperbolicity properties]
The Euler equations for isothermal relativistic fluids
are strictly hyperbolic in the region $|v|< 1/\eps$ for all $\rho \geq 0$ (i.e.~even in 
the presence of vacuum singularities in the mass density field), 
but strict hyperbolicity fails as $v \to \pm 1/\eps$ (i.e.~in 
the presence of light speed singularities in the velocity field). 
\end{lemma} 

In contrast, for polytropic perfect fluids  
the Euler equations also fail to be strictly hyperbolic at the vacuum $\rho=0$.

The $v$-derivatives of the eigenvalues considered as functions of $(\rho, v)$ are 
$$
\lambda_{1v}
= \frac{1- \eps^2 }{ \left( 1 - \eps^2 \, v \right)^2},
\qquad 
\lambda_{2v}
=\frac{1- \eps^2} { \left( 1+ \eps^2 \, v \right)^2},  
$$
while $\lambda_{1 \rho} = \lambda_{2\rho} = 0$. 
Their derivatives along the characteristic fields are
\be
\nonumber 
\begin{aligned} 
&  \nabla \lambda_1 \cdot r_1 = \frac{1 - \eps^2}{1 + \eps^2}\frac1{\rho (1 - \eps^2 v)^2} > 0,
\\
&   \nabla \lambda_2 \cdot r_2 = \frac{1 - \eps^2}{1 + \eps^2} \frac1{\rho (1 + \eps^2 v)^2} > 0, 
\end{aligned}
\ee
while in terms of the Riemann invariants we have  
\be
\label{BP.speed-riemann}
\begin{aligned}  
& \lambda_{1w}=\lambda_{1v}v_w+\lambda_{1\rho}\rho_w
=
  \frac{(1 - \eps^2) \, (1- \eps^2 \, v^2) }{ 2 \, \left( 1- \eps^2 \,  v \right)^2} = \lambda_{1z}, 
\\
& 
\lambda_{2z}=\lambda_{2v}v_z+\lambda_{2\rho}\rho_z
=   \frac{(1 - \eps^2) \, (1 - \eps^2 \, v^2) }{ 2 \, \left( 1+ \eps^2 \, v \right)^2} = \lambda_{2w}.  
\end{aligned} 
\ee

Hence, we conclude:

\begin{lemma}[Genuine nonlinearity property]
The Euler equations for isothermal relativistic fluids admit 
two genuinely nonlinear characteristic fields in the domain 
$\rho \geq 0$, $|v|<1/\eps$. However, the genuine nonlinearity property fails in the limit $v \to \pm 1/\eps$.
\end{lemma}

In contrast, for polytropic perfect fluids the Euler equations also fail to be genuinely nonlinear at the vacuum.

We point out that the wave speeds and Riemann invariants are smooth functions even as $\eps \to 0$: 
\be
\label{DL-eigenv}
\aligned 
& \lam_1 = (v - 1) ( 1 + \eps^2 \, v + O(\eps^4 v^2)), 
\\
& \lam_2 = (v + 1) ( 1 - \eps^2 \, v + O(\eps^4 v^2)),
\endaligned 
\ee
and 
\be
\label{DL-ubar}
\ubar(v) = v \, \left(1 + {\eps^2 \over 3} v^2 + O(\eps^4 v^4) \right).  
\ee
Note that \eqref{DL-eigenv} is uniform in the whole interval $v \in [-1/\eps, 1/\eps]$, while 
the remainder in \eqref{DL-ubar} blows-up when $v$ approaches the light speed. 
  
  
\section{Derivation of additional conservation laws} 

\subsection{Entropy equation} 

From the equations  \eqref{IN.Euler2} and for smooth solutions 
we can derive additional conservation laws, which will play a central role in the existence theory.  
By definition, a pair of mathematical entropy $\Ucal=\Ucal(\rho,v)$ and 
entropy-flux $\Fcal=\Fcal(\rho,v)$ 
provides a conservation law satisfied by all smooth solutions of \eqref{IN.Euler2}. 
The entropy pairs are determined by the compatibility conditions
$$ 
\nabla \Fcal \cdot r_j = \lam_j \, \nabla \Ucal \cdot r_j, \quad j=1,2. 
$$
Expressing $\Ucal, \Fcal$ as functions of $w, z$ and relying on the properties of the Riemann invariants, 
these conditions are equivalent to  
\be
\Fcal_w = \lam_2 \, \Ucal_w, \qquad \Fcal_z = \lam_1 \, \Ucal_z, 
\label{BP.entropypair1}
\ee 
and imply an equation satisfied by the entropy $\Ucal=\Ucal(w,z)$ only: 
$$
\left( \lam_1 \, \Ucal_z \right)_w = \left( \lam_2 \, \Ucal_w \right)_z. 
$$
That is, $\Ucal$ satisfies 
\be 
\Ucal_{wz} + \frac{\lam_{2z} }{ \lam_2 - \lam_1} \, \Ucal_w - \frac{\lam_{1w} }{ \lam_2 - \lam_1} \, \Ucal_z = 0,  
\label{BP.entropy-equation}
\ee 
which we will refer to as the {\sl entropy equation.}

Using the formulas \eqref{BP.speed-riemann}, the coefficients in \eqref{BP.entropy-equation} are found to be 
$$
\aligned 
& \frac{\lambda_{1w}}{\lambda_2-\lambda_1}
=  \frac{(1 - \eps^2) \, (1+\eps^2 \, v) }{ 4 (1-\eps^2 \, v)}, 
\\
& \frac{\lambda_{2z}}{\lambda_2-\lambda_1}
=  \frac{(1 - \eps^2) \, (1 - \eps^2 \, v) }{ 4 (1 + \eps^2 \, v)},
\endaligned 
$$
and thus \eqref{BP.entropy-equation} becomes  
\be
\Ucal_{wz} +  \bbar (w+z) \, \Ucal_w + \abar (w+z) \, \Ucal_z = 0,
\label{BP.equationhyperbolic}
\ee
where the coefficients depend on $w+z$ only. We have set 
$$ 
\begin{aligned}  
& a(v) := -   \frac{(1 - \eps^2 )(1 + \eps^2 v )}{4 (1 - \eps^2 v )},
\\
& b(v) :=   \frac{(1 - \eps^2 )(1 - \eps^2 v )}{4 (1 + \eps^2 v)} = - a(-v), 
\end{aligned} 
$$ 
in which $v = \vbar (w + z)$, and $\abar := a \circ \vbar, \, \  \bbar := b \circ \vbar$, therefore 
\be 
\label{3.5} 
\aligned  
& \abar (\xi ) =  - \frac{1 - \eps^2}4 \, \frac{1 - \eps + (1 + \eps) e^{\eps \xi}}{1 + \eps + (1 - \eps) e^{\eps \xi}},
\\
& \bbar (\xi ) =  \frac{1 - \eps^2}4 \, \frac{1 + \eps + (1 - \eps) e^{\eps \xi}}{1 - \eps + (1 + \eps) e^{\eps \xi}} = - \abar ( - \xi ). 
\endaligned 
\ee

The equation \eqref{BP.equationhyperbolic} is a linear hyperbolic equation
with smooth coefficients, and its solutions are generated by the  
{\sl Riemann function} $\Rcal(w',z'; w,z)$, defined for each fixed $(w,z)$ by the 
{\sl Goursat problem} associated with the adjoint operator: 
\be\label{BP.Goursat1}
\begin{aligned}  
& \Rcal_{w'z'} -\big( \bbar (w'+z') \, \Rcal\big)_{w'} - \big(\abar (w'+z') \, \Rcal\big)_{z'} = 0, 
\\
& \Rcal_{w'} (w',z;w,z) = \abar(w'+z) \, \Rcal(w',z; w,z) \quad \text{ on the line } z'=z, 
\\
& \Rcal_{z'} (w,z';w,z) = \bbar (w+z') \, \Rcal(w,z'; w,z) \quad \text{ on the line } w'=w,
\\
& \Rcal (w,z;w,z)=1. 
\end{aligned} 
\ee 
The Riemann function allows us to solve the general characteristic value problem
\be\label{BP.Goursat3} 
\begin{aligned}  
& \Ucal_{wz} +  \bbar (w+z) \, \Ucal_w + \abar (w+z) \, \Ucal_z = g, 
\\
& \Ucal(w,z') = \varphi(w) \quad \text{ on the line } z=z', 
\\
& \Ucal(w',z) = \psi(z) \quad \text{ on the line } w=w',
\end{aligned}
\ee
where $\varphi, \psi$ are prescribed characteristic data and $g$ is a given source. 
Indeed, we have 
$$
\begin{aligned}  
\Ucal(w,z) := 
&    \frac{1 }{ 2} \, \varphi(w) \, \Rcal(w,z'; w,z)  + \frac{1 }{ 2} \, \psi(z)  \, \Rcal(w',z; w,z) 
\\
&      + \int_{w'}^w  \frac{1 }{ 2} \, \Rcal(w'',z'; w,z) \, \varphi_w(w'') \, dw''
 \\ 
& +    \int_{w'}^w   \Big( (\abar (w''+z') \, \Rcal(w'',z'; w,z) -\frac{1 }{ 2} \, \Rcal_{w'}(w'',z'; w,z) \Big) \,\varphi(w'') \, dw''
\\
&      + \int_{z'}^z  \frac{1 }{ 2} \, \Rcal(w',z''; w,z) \, \psi_z(z'') \, dz''
 \\ 
& +    \int_{z'}^z  \Big( (\bbar (w'+z'') \, \Rcal(w',z''; w,z) - \frac{1 }{ 2} \, \Rcal_{z'}(w',z''; w,z) \Big) \,\psi(z'') \, dz''
\\
& +   \int_{w'}^w   \int_{z'}^z  g(w'',z'') \, \Rcal(w'',z''; w,z)  \, dw'' dz''. 
\end{aligned}  
$$ 


\subsection{Non-relativistic limit}

We assume first that $\eps =0$. Formally, when $\eps \to 0$ in \eqref{BP.equationhyperbolic} we obtain 
$$
\Ucal^0_{wz}+ \frac{1 }{ 4} \, \left( \Ucal_w^0 - \Ucal_z^0 \right) = 0. 
$$
The Riemann function associated with this equation was constructed by LeFloch and Shelukhin (2005): 
\be
\Rcal^0(w',z'; w,z) = e^{((w-w') - (z-z'))/4} \, f^0((w-w') \, (z-z')), 
\label{BP.Riemann-function}
\ee 
where the function $f^0 = f^0(m)$ is related to the Bessel function of order $0$ and can be 
characterized as the solution to the ordinary differential equation
$$ 
\aligned 
& m \, f_{mm}^0 + f_m^0 + f^0/16 = 0, 
\\
& f^0(0) = 1, \quad f^0_m(0) = -1/16. 
\endaligned  
$$
One can check that, for every fixed $(w,z)$, the function $(w',z') \mapsto \Rcal^0(w',z'; w,z)$
defined in \eqref{BP.Riemann-function} is the unique solution of the Goursat problem
$$ 
\begin{aligned}  
&  \Rcal^0_{w'z'} - \frac{1}{4} \, \left( \Rcal^0_{w'} - \Rcal^0_{z'} \right) = 0, 
\\
& \Rcal^0(w',z;w,z) = e^{(w-w')/4} \quad \text{ on the line } z'=z, 
\\
& \Rcal^0(w,z';w,z) = e^{-(z-z')/4} \quad \text{ on the line } w'=w.  
\end{aligned}
$$ 

The following description of the entropy kernel was also established. 

\begin{theorem} [Entropy kernel of the non-relativistic Euler equations]
Consider the isothermal non-relativistic Euler equations 
\eqref{IN.Euler-non} (with $k=1$ after normalization). 
Then the function 
$$
\chi^0(w,z) = \begin{cases}
\Rcal^0(0,0;w,z) = e^{(w - z)/4} \, f^0 (wz),    &   w \, z \leq 0, 
\\
0,             & w \, z > 0, 
\end{cases}
$$
is a fundamental solution of the entropy equation
$$
\chi^0_{wz}+ \frac{1 }{ 4} \, \left( \chi_w^0 - \chi_z^0 \right) = - 2 \delta_{w= 0} \otimes \delta_{z=0}.  
$$
It is a function with bounded variation and  
\be
\aligned 
\chi^0_w (w,z)  
& = e^{-z/4} \,  \big( - (\sgn z) \,\delta_{w=0} + e^{w/4} ( \frac{1 }{ 4} f^0(wz) + z \, f_m^0(wz) ) \, \one_{wz < 0}\big), 
\\
\chi^0_z (w,z) 
& = e^{w/4} \, \big( - (\sgn w) \, \delta_{z = 0} + e^{- z/4} ( - \frac{1 }{ 4} f^0(wz) + w \, f_m^0 (wz) ) \, \one_{wz < 0} \big), 
\endaligned 
\label{BP.singularities1}
\ee
where $\one_{wz< 0}$ denotes the characteristic function of the set $\{wz < 0\}$.
\end{theorem}

It is is natural to define the entropy kernel by
imposing data on the line $\rho=1$. (Observe that $w=z=0$ correspond to $(\rho,v) = (1, 0)$.)  
Hence, the mathematical entropies of \eqref{IN.Euler-non} that vanish on the vacuum are given by the formula
\be 
\Ucal(w,z) = \int_\RR \chi^0(w-s,z-s) \, \psi(s) \, ds, 
\label{BP.allentropies0}
\ee
valid in each region $\rho<1$ and $\rho>1$ (as they avoid the point mass at $w=z=0$), 
where $\psi:\RR \to \RR$ is an arbitrary, integrable function. A similar formula hold for the entropy flux. 
 
We recall that entropy solutions to the non-relativistic Euler
equations admit bounded Riemann invariants, and therefore satisfy 
the restriction $|v \pm \ln \rho | \leq C$ in the physical variables. 
In particular, the argument of the function $f^0$ in the definition of $\chi^0$ remains in a compact set. 
Note also that the entropy kernel satisfies 
$$
|\chi^0(\rho,v) | \lesssim \rho^{1/2}.
$$  
 
 
\subsection{Entropy kernel}
 
We now return to the general model with $\eps \neq 0$. One of the main results 
in the present paper is the following characterization of the entropy kernel.

\begin{theorem}[Entropy kernel of the relativistic Euler equations]
\label{EN-1}
Consider the Euler equations for isothermal relativistic fluids \eqref{BP.Euler3}.   
The function 
$$
\chi(w,z) = \begin{cases}
\Rcal(0,0;w,z),   &   w \, z \leq 0, 
\\
0,             & w \, z > 0, 
\end{cases}
$$
is a fundamental solution of the entropy equation 
$$
\chi_{wz} +  \bbar (w+z) \, \chi_w + \abar (w+z) \, \chi_z = - 2 \delta_{w= 0} \otimes \delta_{z=0}.  
$$
It is solely a function of bounded variation and  
\be\label{BP.singularities10}
\aligned
\lim_{z \to 0, wz <0} \chi(w,z) & = \frac{1+\eps + (1-\eps) e^{\eps w}}2 e^{(1-\eps)^2 w/4},\\
\lim_{w \to 0, wz<0} \chi(w,z) & = \frac{1-\eps + (1+\eps) e^{\eps z}}2 e^{-(1+\eps)^2 z/4},
\endaligned
\ee 
\be
\aligned 
\chi_w (w,z)  & = -  (\sgn z) \frac{1 - \eps + (1 + \eps )e^{\eps z}}2 e^{ -  (1 + \eps)^2 z/4}  
\delta_{w=0} + C_1(w,z) \, \one_{wz < 0}, 
\\
\chi_z (w,z) & = - (\sgn w) \frac{1 + \eps + (1 - \eps ) e^{ \eps w}}2 \, e^{ (1 - \eps)^2 w/4} \, 
\delta_{z = 0} + C_2(w,z) \, \one_{wz < 0}, 
\endaligned
\label{BP.singularities11}
\ee
where $C_1, C_2$ are smooth functions.
In the physical variables $\rho,v$, the kernel $\chi=\chi(\rho,v)$ has compact support in the variable 
$v \in (-1/\eps,1/\eps)$
(for every fixed value $\rho$), and satisfies 
$$
\begin{array}{ll}   
& \lim\limits_{\rho \to 1} \chi(\rho,\cdot) = 0, \qquad  
       \lim\limits_{\rho \to 1\pm} \chi_\rho(\rho,\cdot) = \pm \frac1{1+\eps^2} \delta_{v=0} 
       \quad \text{in the weak sense,}
\\
&   |\chi(\rho,u)| \lesssim \rho^\alpha, \qquad 
\alpha := \frac{(1+\eps)^2 }{ 2 (1 + \eps^2)}.     
\end{array} 
$$ 
\end{theorem}

In the course of the proof of this theorem we will also show
\be
\label{BP.singularities12}
\aligned
C_1(w,0)  & = \frac{1 - \eps^2}{4} \frac{1 - \eps + (1 + \eps) e^{\eps w}}2 \, e^{ (1 - \eps)^2 w / 4},
\\ 
C_2(0,z) & = - \frac{1 - \eps^2}{4} \frac{1 + \eps + (1 - \eps) e^{\eps z}}2 \, e^{ -(1 + \eps)^2 z/ 4},
\endaligned 
\ee 
and
\be\label{BP.singularities13}
\aligned 
& C_1(0,z) = {1 - \eps^2 \over 4} \left( 1 - \frac{1-\eps^2}4 z \right) \, 
\frac{1 - \eps + (1 + \eps )e^{\eps z}}2 e^{ -  (1 + \eps)^2 z/4} , 
\\
& C_2(w,0) = - {1 - \eps^2 \over 4} \left( 1 + \frac{1-\eps^2}4 w \right)\frac{1 + \eps + (1 - \eps ) e^{ \eps w}}2 \, e^{ (1 - \eps)^2 w/4}. 
\endaligned 
\ee

Observe that the kernel $\chi$ is only H\"older continuous at the vacuum; the H\"older 
exponent converges to $1/2$ as $\eps \to 0$, allowing us to recover the behavior in $\sqrt{\rho}$
known in the non-relativistic regime. 
(In fact, as $\eps \to 0$,  the expansion of the entropy kernel converges term by term 
to the one of the non-relativistic case.) 
Next, we introduce:

\begin{definition}[Notion of weak entropy]
A pair of continuous maps $(\Ucal,\Fcal) : \RR_+ \times (-1/\eps, 1/\eps) \to \RR^2$ is called a 
{\rm weak entropy pair}
of the isothermal relativistic Euler equations if the partial differential equations \eqref{BP.entropypair1}
hold in the sense of distributions and, moreover, $\Ucal$ vanishes on the vacuum within any tame region, 
in the following sense: for any $M>0$ and $\eps \in (0,1)$ 
there exist some constant $C>0$ such that 
\be
|\Ucal(\rho,v)| + |\Fcal(\rho,v)| \leq C \, \rho^\alpha,  \qquad (\rho,v) \in \Tcal_\eps(M).  
\label{AS.vanish}
\ee
\end{definition}

For instance, the conservative and flux variables in \eqref{IN.Euler2} do qualify as weak entropies. 
From Theorem~\ref{EN-1} we deduce: 

\begin{corollary}[Entropy pairs of the relativistic Euler equations] 
In each of the regions $\rho<1$ and $\rho>1$ the formula 
\be
\Ucal(w,z) = \int_\RR \chi(w-s,z-s) \, \psi(s) \, ds
\label{BP.allentropies1}
\ee
determines a family of weak entropies, where $\psi:\RR \to \RR$ is an arbitrary integrable function. 
\end{corollary}


\subsection{Entropy flux}

Given an entropy $\Ucal$ (with sufficient decay) we deduce from \eqref{BP.entropypair1}, 
that the associated entropy flux is 
$$   
\aligned 
\Fcal(w,z) 
& =  \int_{-\infty}^w \lambda_2(w',z) \, \Ucal_w(w',z) \, dw'
\\
& = \lambda_2(w,z) \, \Ucal(w,z) -  \int_{-\infty}^w \lambda_{2w}(w',z) \, \Ucal(w',z) \, dw'. 
\endaligned 
$$ 
A similar formula hold with $\lambda_2$ replaced by $\lambda_1$ and, therefore, 
$$
\begin{aligned}  
\Fcal(w,z) =  &   \frac{1 }{ 2} \left(\lambda_1(w,z) + \lambda_2(w,z)\right) \, \Ucal(w,z)
 \\
&   - \frac{1 }{ 2} \int_{-\infty}^w \lambda_{2w}(w',z) \, \Ucal(w',z) \, dw'
  + \frac{1 }{ 2} \int_z^\infty \lambda_{1z}(w,z') \, \Ucal(w,z') \, dz'. 
\end{aligned}  
$$
These identities allow us to define the entropy flux kernel, denoted below by $\sigma=\sigma(\rho,v,s)$ 
from the entropy kernel $\chi$,
and in turn to compute the entropy flux via the formula 
\be
\label{forml}
\Fcal (\rho, v) = \Fcal (w,z) = \int_\RR \sigma (w - s, z - s) \,  \psi (s) \, ds.
\ee

More precisely, the condition \eqref{BP.entropypair1} reads 
\be
\label{BP.entropypair2}
\sigma_w = \lambda_2 \chi_w, \qquad \sigma_z = \lambda_1 \chi_z.
\ee
By taking advantage of the Lorentz invariance property of the relativistic Euler equations
we can decompose the entropy flux as follows: 
\be
\label{lorentz1} 
\sigma(\rho,v,s) 
=: \ubar(v) \, \chi(\rho, v,s) + \sigmas(\rho,v,s),
\ee
or, equivalently, $\sigma(\rho,v,s) = \frac{w +z}2 \, \chi(\rho, v,s) + \sigmas(\rho,v,s)$. 
The key property (which is easily checked from the equations defining the entropy and entropy flux kernel)
is that the dependence of $\chi$ and $\sigma$ with respect to $s$ can be suppressed, since 
$$
\aligned
& \chi(\rho,v,s) = \chi(w-s,z-s),
\\
& \sigmas (\rho,v,s) = \sigmas (w-s,z-s).
\endaligned
$$
Clearly, the condition \eqref{BP.entropypair2} is equivalent to  
\be
\label{lorentz3}
\aligned 
& \sigmas_w = - \frac{w +z}2 \, \chi_w - \frac12 \chi + \lam_2 \, \chi_w, 
\\
& \sigmas_z = - \frac{w +z}2 \, \chi_z - \frac12  \chi + \lam_1 \, \chi_z
\endaligned 
\ee

Moreover we impose the following boundary condition on $\sigmas=\sigmas(w,z)$:
$$
\sigmas(0,0) = 0. 
$$
In particular, the traces along the boundaries $w=0$ and $z=0$ can be determined explicitly from 
the boundary conditions for the Riemann function \eqref{BP.Goursat1} and the expansion 
\eqref{BP.singularities11}. By plugging \eqref{speeds} and \eqref{BP.singularities11},
we obtain the following corollary of Theorem~\ref{EN-1}.

\begin{lemma}[Entropy flux kernel] 
\label{corolla} 
\be
\label{BP.entropypair3} 
\aligned 
 \sigmas_w = & C_5(z) \, e^{-(1 + \eps)^2 z/4} \, \delta_{w=0}
  + C_3(w,z) \one_{wz < 0}, 
\\ 
\sigmas_z = & C_6(w) \, e^{(1 - \eps)^2 w/4} \, \delta_{z = 0}  + C_4(w,z) \one_{wz < 0},
\endaligned 
\ee
where 
$$\aligned
C_3 (w,z) & = \big( \lambda_2 (w,z) - \frac{w+z}2 \big) \, C_1 (w,z) - \frac12 \chi (w,z),
\\
C_4 (w,z) & = \big( \lambda_1 (w,z) - \frac{w+z}2 \big) \, C_2 (w,z) - \frac12 \chi (w,z),
\endaligned
$$
and 
$$
\aligned 
& C_5(z) := |z| \, \frac{1 - \eps + (1 + \eps ) e^{\eps z}}4 +  (\sgn z) \, \frac{1 - \eps - (1 + \eps) e^{\eps z}}{2 \eps}, 
\\
& C_6(w) := |w| \, \frac{1 + \eps + (1 - \eps) e^{\eps w}}4 +  (\sgn w) \, 
\frac{1 + \eps - (1 - \eps) e^{\eps w}}{2 \eps}.
\endaligned 
$$
\end{lemma}

Finally, it is not difficult to check that, again in each of the regions $\rho<1$ and $\rho>1$, the formula 
\be
\aligned
\Fcal(w,z) 
& = \int_\RR \sigma (w-s,z-s) \psi (s)  \, ds\\
& = \int_\RR \big( \ubar(v) \chi(w-s,z-s) + \sigma^\sharp(w-s,z-s) \big) \, \psi(s) \, ds
\endaligned
\label{BP.allentropies3}
\ee
determines the entropy flux associated with \eqref{BP.allentropies1}. 


\section{Technical estimates}

In this section we give a proof of Theorem~\ref{EN-1}. 
The existence of the Riemann function is standard, as it is determined by a Goursat problem
for a linear hyperbolic equation with regular coefficients. The main issues to be dealt with are 
the behavior of $\chi$ near $\rho=0$ and near ${v = \pm 1/\eps}$, as well as the expansion 
of $\chi$ along the boundary of its support.

\begin{lemma} The traces of the kernel along its support ($wz < 0$) are given by 
\be
\label{chi.value0} 
\chi (w,z) \to
\begin{cases}
A(w):=\frac{1 + \eps + (1 - \eps) e^{\eps w}}2 \, e^{ (1 - \eps)^2 w / 4}, & z \to 0,
\\
B(z):=\frac{1 - \eps + (1 + \eps )e^{\eps z}}2 \, e^{ - (1 + \eps)^2 z / 4},  & w \to 0, 
\end{cases}
\ee
and, moreover, within any compact set in $\rho$ 
$$
\chi(\rho,v) \lesssim \rho^\alpha. 
$$ 
\end{lemma}

\begin{proof} We need to integrate out the boundary conditions arising in \eqref{BP.Goursat1}. 
The second differential equation in \eqref{BP.Goursat1} implies 
$$
\ln \vert \Rcal (w',z;w,z)  \vert + C_0 
= - \frac{1 - \eps^2}4 \Big( - (w' + z)
+ 2 \int^{w'} \frac{1 }{ 1 - \eps^2 v} dw' \Big)
$$
for some constant $C_0$, where
$$
\aligned 
\int^{w'} \frac{1 }{ 1 - \eps^2 v} dw' = &   \int^{w'} \frac{1 }{ 1 - \eps \tanh \left( \eps (w' + z)/2\right)} dw' 
\\
= &    \frac{2}{1 - \eps^2} \ln \Big( (1 - \eps ) e^{\eps (w' + z)}
+ 1 + \eps \Big) + \frac1{1 + \eps}  (w' + z).
\endaligned
$$
The third equation in \eqref{BP.Goursat1} implies
$$
\ln \vert \Rcal (w,z';w,z) \vert + C_0 = \frac{1 - \eps^2}4 \Big( - (w + z')
+ 2 \int^{z'} \frac{1 }{ 1 + \eps^2 v} dz' \Big)
$$
 for some constant $C_0$, where
$$
\aligned 
\int^{z'} \frac{1 }{ 1 + \eps^2 v} dz'
= &   \int^{z'} \frac{1 }{ 1 + \eps \tanh \left( \eps (w + z')/2\right)} dz' \\
= &   - \frac{2}{1 - \eps^2} \ln \Big( (1 + \eps ) e^{\eps (w + z')} + 1 - \eps \Big) + \frac1{1 - \eps} (w + z').
\endaligned 
$$

Since $\Rcal (w,z;w,z)=1$, we find 
$$
\aligned 
& \Rcal (w',z';w,z) 
\\
& = 
\begin{cases}
e^{- (1 - \eps)^2 (w' - w) / 4}  \, 
\frac{1 - \eps + (1 + \eps )e^{ - \eps (w + z)}}{(1 - \eps)e^{\eps (w' - w)} + (1 + \eps )e^{-\eps (w + z)}}
=: A(w'; w,z),     &  z' = z,
\\
e^{ (1 + \eps)^2 (z' - z) / 4}  \, \frac{1 + \eps + (1 - \eps )e^{-\eps (w + z)}}{(1 + \eps)e^{\eps (z' - z)} + (1 - \eps )e^{-\eps (w + z)}}=: B(z';w,z), & w' = w, 
\end{cases}
\endaligned 
$$
and in particular 
$$
\Rcal (0,0;w,z) = 
\begin{cases}
e^{ (1 - \eps)^2 w / 4}  \, \frac{1 + \eps + (1 - \eps )e^{\eps w}}{2}=A(0;w,0)= A(w),     &  z = 0,
\\
e^{ - (1 + \eps)^2 z / 4}  \, \frac{1 - \eps + (1 + \eps )e^{\eps z}}{2}=B(0;0,z)= B(z),
 & w = 0,
\end{cases}
$$
which provides the desired behavior on $\chi$ along its support.

Next, we write along the boundary $z=0$ (corresponding to $u=R$ and $w=2R$) 
$$
A(w) = {1 \over 2} \rho^{(1-\eps)^2 \over 2(1 + \eps^2)} \, 
\big( 1+ \eps + (1-\eps) \rho^{{2\eps \over 1 + \eps^2}}\big), 
$$
which shows that $\chi(\rho,v) \lesssim \rho^\alpha$. On the boundary $w=0$ we have $u=-R$ and $z=-2R$ and we find 
a similar estimate $\chi(\rho,v)  \lesssim \rho^\alpha$ since  
$$
B(z) = {1 \over 2} \rho^{(1+\eps)^2 \over 2(1 + \eps^2)} \, \big( 1- \eps + (1+\eps) \rho^{-{2\eps \over 1 + \eps^2}}\big).
$$
Since $\chi$ is smooth, the estimate remains valid in any compact region. 
\end{proof}

\begin{lemma} 
The entropy kernel $\chi (w,z)=\Rcal (0,0;w,z)\one_{wz < 0}$ is a fundamental solution of the entropy 
equation 
$$
\chi_{wz} + \bbar (w+z) \chi_w + \abar (w+z) \chi_z = -2 \delta_{w=0} \otimes \delta_{z=0}.
$$
The traces of the derivatives $\chi_w$ and $\chi_z$ along the boundaries $z=0$ and $w=0$ (while keeping $wz<0$) are
\be
\label{chi.value1a}
\aligned
\lim_{z \to 0}\chi_w(w,z) &=
C_1(w,0) = A_w(w)
\\
& = \frac{1 - \eps^2}{4} \frac{1 - \eps + (1 + \eps) e^{\eps w}}2 \, e^{ (1 - \eps)^2 w / 4} 
  = - \abar (w) A(w),
\endaligned 
\ee
and
\be
\label{chi.value1b}
\aligned
\lim_{w \to 0}\chi_z(w,z) &=
C_2(0,z) = B_z(z)
\\ 
& = - \displaystyle\frac{1 - \eps^2}{4} \frac{1 + \eps + (1 - \eps) e^{\eps z}}2 \, e^{ -(1 + \eps)^2 z/ 4}
  = - \bbar (z) B(z). 
\endaligned 
\ee 
\end{lemma}

\begin{proof} In view of \eqref{chi.value0}, we obtain 
\be
\label{BP.singularities2}
\aligned 
 \chi_w(w,z) & = \Rcal_w (0,0;w,z) \one_{wz < 0} - (\sgn z) B(z) \delta_{w=0}
\\
& = - (\sgn z) \, \frac{1 - \eps + (1 + \eps) e^{\eps z}}2 \,  e^{-(1+\eps)^2 z/4}    \delta_{w=0} + C_1(w,z) \one_{wz < 0},\\
 \chi_z(w,z) & = \Rcal_z (0,0;w,z) \, \one_{wz < 0} - (\sgn w) \, A(w) \, \delta_{z=0}
 \\
& = - (\sgn w) \, \frac{1 + \eps + (1 - \eps) \, e^{\eps w}}2 \, e^{(1-\eps)^2 w/4}  \delta_{z=0} + C_2(w,z) \, \one_{wz < 0}. 
\endaligned 
\ee
On the other hand, from \eqref{chi.value0} we can determine  
\be
\nonumber
\aligned
C_1(w,0) & = A_w(w)
\\
& = \frac{d}{dw} \Big(\frac{1 + \eps + (1 - \eps) e^{\eps w}}2 \, e^{ (1 - \eps)^2 w / 4}\Big)
\\
& = (1 - \eps^2) \frac{1 - \eps + (1 + \eps) e^{\eps w}}8 \, e^{ (1 - \eps)^2 w / 4} 
  = - \abar (w) A(w)
\endaligned 
\ee
and 
\be
\nonumber
\aligned
C_2(0,z) & = B_z(z)\\ 
& = \frac{d}{dz} \left(\frac{1 - \eps + (1 + \eps )e^{\eps z}}2 \, e^{ - (1 + \eps)^2 z / 4}\right) 
\\
& = - (1 - \eps^2) \frac{1 + \eps + (1 - \eps) e^{\eps z}}8 \, e^{ -(1 + \eps)^2 z/ 4}
  = - \bbar (z) B(z). 
\endaligned 
\ee

We differentiate \eqref{BP.singularities2} once again, and in view of \eqref{chi.value1a} and \eqref{chi.value1b} we obtain
\begin{align*} 
\chi_{wz}(w,z)
& = \Rcal_{wz} (0,0;w,z) \, \one_{wz < 0}  - (\sgn w) \, \Rcal_w(0,0;w,0) \, \delta_{z=0}
\\ 
& \quad  - (\sgn z) \, B_z(z) \, \delta_{w=0} - 2 \, B(0) \, \delta_{w=0} \otimes\delta_{z=0}
\\
& = \Rcal_{wz} (0,0;w,z) \, \one_{wz < 0}  - (\sgn w) \, A_w(w) \, \delta_{z=0}
\\
& \quad - (\sgn z) \, B_z (z) \, \delta_{w=0}
- 2  \delta_{w=0}\otimes\delta_{z=0}
\\
& =  \Rcal_{wz} (0,0;w,z) \, \one_{wz < 0}  + (\sgn w) \, \abar (w) A(w) \, \delta_{z=0}
\\
& \quad + (\sgn z) \, \bbar (z) \, B(z) \, \delta_{w=0} - 2  \, \delta_{w=0}\otimes\delta_{z=0}, 
\end{align*}
which allows us to compute 
\begin{align*}
& \chi_{wz} + \bbar (w+z) \, \chi_w + \abar (w+z) \, \chi_z
\\
& = \Rcal_{wz} (0,0;w,z) \, \one_{wz < 0} + (\sgn w) \, \abar (w) \, A(w) \, \delta_{z=0}
\\
& \quad + (\sgn z) \, \bbar (z) \, B(z) \, \delta_{w=0} - 2 \, \delta_{w=0}\otimes\delta_{z=0}
\\
& \quad + \bbar (w+z) \, \left(\Rcal_w (0,0;w,z) \, \one_{wz < 0} - (\sgn z) \, B(z)
\delta_{w=0}\right)
\\
& \quad + \abar (w+z) \, \left(\Rcal_z (0,0;w,z) \, \one_{wz < 0} - (\sgn w) A(w) \, \delta_{z=0}\right)
\\
& = \left( \Rcal_{wz} (0,0;w,z) + \bbar (w+z) \, \Rcal_w(0,0;w,z)
      + \abar (w+z) \, \Rcal_z (0,0;w,z)\right) \, \one_{wz < 0} 
\\
& \quad + (\sgn w) \, \left( \abar (w) - \abar (w+z) \right) \, A(w) \, \delta_{z=0}
\\
& \quad + (\sgn z) \, \left( \bbar (z) - \bbar (w+z) \right) \, B(z) \, \delta_{w=0}
   -2 \, \delta_{w=0}\otimes\delta_{z=0}
\\
& = -2  \delta_{w=0}\otimes\delta_{z=0}. 
\end{align*} 
\end{proof}

\begin{lemma} 
The traces of the derivatives $\chi_w$ and $\chi_z$ along the boundaries $w=0$ and 
$z=0$ (while keeping $wz<0$)  are  
\be
\label{BP.singularities3a}
\aligned
\lim_{w \to 0}\chi_w(w,z) & = 
 C_1(0,z) 
 \\
  & = \frac{1 - \eps^2}4 \left( 1 - (1-\eps^2) z/4 \right) \frac{1 - \eps + (1 + \eps )e^{\eps z}}2 e^{ -  (1 + \eps)^2 z / 4}\\
& = {1 - \eps^2 \over 4}\left( 1 - (1-\eps^2) z/4 \right) \chi(0,z)
\endaligned 
\ee
and 
\be
\label{BP.singularities3b}
\aligned 
\lim_{z \to 0}\chi_z(w,z)& =
C_2(w,0) 
\\ & = - {1 - \eps^2 \over 4} \left( 1 + (1-\eps^2) w/4 \right) \frac{1 + \eps + (1 - \eps ) e^{ \eps w}}2 \, e^{ (1 - \eps)^2 w / 4}\\
& =  - {1 - \eps^2 \over 4}\left( 1 + (1-\eps^2) w/4 \right) \chi(w,0). 
\endaligned 
\ee
\end{lemma}

\begin{proof} 
We give the proof for the boundary $z=0$ with $w>0$, the calculation for the other boundary being similar.
Using \eqref{BP.Goursat3} and the expressions \eqref{3.5} we write 
$$
{d \over dw} \chi_z(w,0-)  + \abar (w) \, \chi_z(w,0-) = - \bbar (w) \, \chi_w(w,0-),  
$$
where the coefficients $a$ and $b$ can be regarded as functions of $w+z=w$ and the right-hand side is already known.  
By integrating the above equation we obtain 
\be
\label{laformule}
\aligned
& \chi_z(w,0-) 
\\
& = -\int_0^w \bbar(w') \, \chi_w(w',0) \, e^{\int_w^{w'} \abar(w'') \, dw''} dw'
    + \chi_z(0,0) \, e^{-\int_0^w \abar(w') \, dw'},  
\endaligned
\ee
where  $\chi_z(0,0)$ stands for $\lim_{w \to 0+} \chi_z(w,0-)$.

In view of \eqref{chi.value1b}, $\chi_z(0,0) = -(1 - \eps^2)/4$. In view of \eqref{chi.value1a} and $A(0)=1$,
$$
e^{- \int_0^w \abar(w') \, dw' } = A(w) = \frac{1 + \eps + (1 - \eps) e^{\eps w}}2 \, e^{ (1 - \eps)^2 w / 4} 
$$
and 
$$
e^{\int_w^{w'} \abar(w'') \, dw''} = \frac{A(w)}{A(w')} = {1 + \eps + (1 - \eps) e^{\eps w} \over 1 + \eps + (1 - \eps) e^{\eps w'}} \, e^{ -(1 - \eps)^2 (w'-w) / 4}. 
$$
Returning to \eqref{laformule}, we find
$$
\aligned
\bbar(w') \, \chi_w(w',0) \, e^{\int_w^{w'} \abar(w'') \, dw''} 
&=-\bbar (w')\, \abar (w') A(w') \,  
\frac{A(w)}{A(w')}\\
&=  {(1 - \eps^2)^2 \over 32} \big(1 + \eps + (1 - \eps) e^{\eps w} \big) \, e^{ (1 - \eps)^2 w / 4}
\endaligned
$$
and 
$$\chi_z (w,0-) = - \frac{1-\eps^2}4 \left( 1 + (1-\eps^2) w/4 \right) A(w),
$$
thus
$$
C_2(w,0)  = - {1 - \eps^2 \over 4} \left( 1 + (1-\eps^2) w/4 \right) \frac{1 + \eps + (1 - \eps ) e^{ \eps w}}2 \, e^{ (1 - \eps)^2 w / 4}. 
$$
\end{proof}

The singularities \eqref{BP.singularities11} of the derivatives of $\chi$ with its values at the boundary \eqref{BP.singularities12}, \eqref{BP.singularities13} follow then from
\eqref{BP.singularities2} and from 
\eqref{chi.value1a}, \eqref{chi.value1b}, \eqref{BP.singularities3a}, and \eqref{BP.singularities3b}. 
Moreover substituting \eqref{BP.singularities11} into
$$
\chi_\rho = (\partial_w - \partial_z) \chi \, R_\rho,
$$
we obtain:

\begin{lemma}
The traces of the derivatives $\chi_\rho$ along $\rho = 1 \pm$ are  
$$
\lim_{\rho \to 1 \pm} \chi_\rho (\rho ,v) = \pm \frac1{1 + \eps^2} \delta_{v=0}.
$$ 
\end{lemma}

\begin{lemma}[Derivatives of the entropy kernel and the entropy-flux kernel]\label{BP.singularities4}
We have  
\be\label{BP.singularities4a}
\aligned 
& \chi_s (\rho,v,s)
\\
&   = X^{(1)}(w-z) \left( - \delta_{s=w} +   \delta_{s=z} \right)
 + X^{(2)} (s-w, s-z) \,  \one_{(s-w)(s-z)< 0}, 
\endaligned 
\ee
where 
$$
\aligned 
X^{(1)} (w-z) & :=     \sgn (w - z) \, \frac{1 + \eps + (1 - \eps) \, e^{\eps (w - z)}}2 \, e^{(1 - \eps)^2 \, (w - z)/4}\\
& = \sgn (w - z) \lim_{\eta \to 0, \xi \eta<0} \chi (\xi,\eta)|_{\xi = w-z},
\endaligned 
$$
and 
$$
\aligned 
X^{(2)} (\rho,v,s) & = X^{(2)} (s-w,s-z)  := - \left( C_1 (w-s, z-s) + C_2 (w-s, z-s) \right). 
\endaligned 
$$ 
We have also 
\be
\label{BP.singularities4b} 
\aligned 
& \sigmas_s (\rho,v,s) 
\\
&  = X^{(3)}(w-z) \left( \delta_{s=w} +  \delta_{s = z}\right) 
 + X^{(4)}(s-w, s-z) \one_{(s - w)(s - z) < 0},  
\endaligned 
\ee
where  
$$
\aligned 
X^{(3)} (w-z) & :=   - \left(  \frac{|w - z|}2 \, \frac{1 + \eps + (1 - \eps) e^{\eps( w - z)}}2 \right.\\
   &  \quad  \left.  +   \sgn (w - z)  \, \frac{1 + \eps - (1 - \eps) e^{ \eps (w - z)}}{2\eps} \right)
    \, e^{ (1 - \eps)^2 (w - z)/4}\\
    & = \sgn (w - z) \left( \lambda_1 (w-z,0) 
- \frac{w-z}2 \right) \lim_{\eta \to 0, \xi\eta<0} 
\chi(\xi,\eta)|_{\xi= w-z},
\endaligned 
$$
and 
$$ 
\aligned 
X^{(4)} (\rho,v,s) & = X^{(4)} (s-w,s-z)  :=  - \left( C_3 ( w-s, z-s) +  C_4 (w-s, z-s) \right). 
\endaligned 
$$
\end{lemma}

\begin{proof}
In view of \eqref{BP.singularities11} and \eqref{BP.entropypair3}, we obtain 
$$
\aligned
\chi_s (\rho,v,s) 
= & \frac{\partial}{\partial s} \left( \chi (w-s, z-s) \right) = - (\chi_w + \chi_z) (w-s, z-s)\\
= & \sgn (z - s) \, \frac{1 - \eps + (1 + \eps)e^{\eps (z - s)}}2 e^{- (1 + \eps)^2 (z - s)/4} \, \delta_{w=s}\\  
& - C_1 (w-s, z-s) \, \one_{(w-s)(z-s) < 0}\\
& +   \sgn (w - s) \, \frac{1 + \eps + (1 - \eps) e^{\eps (w - s)}}2 e^{(1 - \eps)^2 (w - s)/4} \, \delta_{z=s}\\
& - C_2(w-s, z-s) \, \one_{(w-s)(z-s) < 0}
\endaligned$$
and thus 
$$
\aligned
\chi_s (\rho,v,s) 
= &  \sgn (w - z) \frac{1 + \eps + (1 - \eps) e^{\eps (w - z)}}2 e^{(1 - \eps)^2 (w - z)/4} \left(- \delta_{s=w} + \delta_{s=z}\right)\\
& - \left( C_1 (w-s,z-s) + C_2 (w-s,z-s) \right) \one_{(w-s)(z-s)<0}.
\endaligned
$$
Defining 
$X^{(1)} (w-z)$
and $X^{(2)} (s-w,s-z)$ as in the lemma 
we obtain the derivative of the entropy kernel \eqref{BP.singularities4a}. 
In view of \eqref{BP.singularities10}, we have also
$$
X^{(1)} (w-z)  = \sgn (w - z) \lim_{\eta \to 0, \xi \eta<0} \chi (\xi,\eta)|_{\xi = w-z}.
$$

Similarly, for the entropy flux kernel we use using \eqref{BP.entropypair3} and obtain
$$
\aligned
& \sigmas_s (\rho,v,s) = \frac{\partial}{\partial s} \left( \sigmas (w-s, z-s)\right)
       = - ( \sigmas_w + \sigmas_z)(w-s, z-s)\\
= & \left( - |z - s| \, \frac{1 - \eps + (1 + \eps) e^{\eps( z - s)}}4  - \sgn (z - s) \, \frac{1 - \eps - (1 + \eps) e^{\eps (z - s)}}{2 \eps}\right)\\
& \, e^{- (1 + \eps)^2 (z - s)/4}\, \delta_{s=w} - 
C_3(w-s, z-s) \one_{(w - s)(z - s) < 0}\\
& + \left( - |w - s| \, \frac{1 + \eps + (1 - \eps) e^{\eps (w - s)}}4 - \sgn (w - s)  \frac{1 + \eps - (1 - \eps) e^{ \eps ( w - s)}}{2 \eps}\right)\\
& \, e^{ (1 - \eps)^2 (w - s)/4} \, \delta_{s = z}
 - C_4(w-s, z-s) \one_{(w - s)(z - s) < 0}\\
= & - \left(  \frac{|w-z|}2 \frac{1 + \eps + (1 - \eps ) e^{\eps (w-z)}}2 + \sgn (w-z) \frac{1 + \eps - (1 - \eps) e^{\eps (w-z)}}{2\eps}\right)\\
& \quad \, e^{(1 - \eps)^2 (w-z)/4} \left( \delta_{s=w} + \delta_{s=z}\right)\\
& \quad -  \left( C_3 (w-s,z-s) + C_4 (w-s,z-s)\right)  \one_{(w-s)(z-s)<0}.
\endaligned
$$
Defining $X^{(3)} (w-z)$ and $X^{(4)} (s-w,s-z)$ as in the lemma
we obtain the derivative of the entropy-flux kernel \eqref{BP.singularities4b}. 
Furthermore, in view of \eqref{speeds}, \eqref{BP.singularities10}
$$\aligned
X^{(3)} & = \sgn (w - z) \big( \lambda_1 (w-z,0) 
- \frac{w-z}2 \big) \lim_{\eta \to 0, \xi\eta<0} 
\chi(\xi,\eta)|_{\xi= w-z}
\\
& = \sgn (-w+z) \big( \lambda_2 (-w+z,0)  - \frac{-w+z}2 \big) \frac{1+\eps + (1-\eps)e^{\eps (w-z)}}2e^{(1-\eps)^2(w-z)/4}. 
\endaligned
$$
\end{proof}


\section{Uniform estimates for the Lax-Friedrichs scheme}
\label{AS-0} 
 
In this section, we follow DiPerna (1983) who considered non-relativistic polytropic fluids, and 
we apply the Lax-Friedrichs scheme. As observed by Hsu, Lin, and Makino (2004), 
DiPerna's arguments carry over to relativistic fluids. Our main purpose is to derive new uniform bounds, 
and establish that Lax-Friedrichs approximations remain in a tame region,  
which is uniquely determined from the initial data of the initial value problem under consideration. 

We begin with: 

\begin{lemma}[Riemann problem]
\label{511} 
Consider the relativistic Euler equations \eqref{IN.Euler2} for $\eps \in (0,1]$.
The Riemann problem corresponding to an initial data made of 
a single jump discontinuity (at $x=0$) separating two constant states $(\rho_l,v_l)$ and 
$(\rho_r,v_r)$ admits a unique self-similar solution  $(\rho,v) = (\rho,v)(x/t)$ satisfying all entropy inequalities.  
The solution of the Riemann problem satisfies the uniform $L^\infty$ bounds in the Riemann invariant variables, 
\be
\label{max-prin}
\aligned 
& 0 \leq W(x/t) \leq W_0:=\sup(W_l, W_r), 
\\
& 0 \leq Z(x/t) \leq Z_0:=\sup(Z_l, Z_r),
\endaligned 
\ee
where the notation $W,Z$ stands for the modified Riemann invariants associated with the variables $(\rho,v)$.
As a consequence, there exists a constant $M>0$ depending only on $(\rho_l,v_l)$ and 
$(\rho_r,v_r)$ such that 
\be
\label{52}
\aligned 
& 0 \leq \rho(x/t) \leq M, 
\\
& 1 - \eps |v(x/t)| \geq (\rho(x/t)/M)^{\eps'}, 
\endaligned 
\ee
where $\eps'$ was defined in the introduction. 
\end{lemma}

\begin{proof} When $\eps \in (0,1)$ the Riemann problem was solved by Smoller and Temple (1993). 
In particular, they established that the shock curves lie in the interior of the rectangular regions 
limited by the integral curves of the characteristic fields. 
This property implies directly the maximum principle \eqref{max-prin}. 
Recall that the rarefaction curves are given by
$$
\aligned
\Rbf_1: & \, u + R = u_l + R_l, \qquad \rho \leq \rho_l,
\\
\Rbf_2: & \, u - R = u_l - R_l, \qquad \rho \geq \rho_l.
\endaligned
$$
Hence on $\Rbf_1$, $w$ is constant and $z$ increases and on $\Rbf_2$, $z$ is constant and $w$ increases.
The shock curves are given by ($i=1,2$) 
$$
\Sbf_i : \, \rho = \rho (v; \rho_l, v_l) 
= \rho_l + \rho_l \, \beta(v, v_l) \, \left( 1 + (-1)^i \sqrt{1 + 2/\beta} \right), 
$$
where
$$
\beta (v, v_l) := \frac{(1 + \eps^2)^2}2 \frac{(v - v_l)^2}{(1 - \eps^2 v^2)(1 - \eps^2 v_l^2)}.
$$

From the definition \eqref{modifiedRinvariants}, the condition \eqref{max-prin} can be rewritten in the form 
$$
\frac{1}{Z_0} \, \rho(x/t)^{1/(1+\eps^2)} \leq \Big(\frac{1 + \eps v(x/t)}{ 1 - \eps v(x/t)} \Big)^{1/(2\eps)}
\leq W_0 \, \rho(x/t)^{-1/(1+\eps^2)}, 
$$
which leads to \eqref{52} with a constant 
$M_\eps:=\left(\sup (W_0, Z_0)\right)^{1 
+ \eps^2}$ depending upon $\eps$.  In fact, from $\rho = (WZ)^{(1 + \eps^2)/2}$, 
$$\rho(x/t)\leq (W_0Z_0)^{(1 + \eps^2)/2}
\leq M_\eps.$$ 
On one hand, $\frac{1 + \eps v(x/t)}{1 
- \eps v(x/t)} \leq \rho(x/t)^{- 2 
\eps/(1 + \eps^2)} W_0^{2 \eps}$ leads to 
$$\aligned
 1 - \eps v(x/t)  & \geq \frac2{1 
+\rho(x/t)^{- 2 \eps/(1 + \eps^2)} 
W_0^{2 \eps}} \geq \frac2{1 + 
\left(\rho(x/t)/ M_\eps\right)^{- 2 
\eps/(1 + \eps^2)}}\\
& \geq \left( \rho(x/t) / M_\eps \right)^{\eps'}
\endaligned
$$
On the other hand, the inequality 
$$
\frac{1 + \eps v(x/t)}{1 - \eps v(x/t)} 
\geq \rho(x/t)^{2\eps/(1 + \eps^2)} Z_0^{-2 \eps}
$$
leads us to
$$\aligned
 1 + \eps v(x/t) & \geq 
\frac2{1+\rho(x/t)^{- 2\eps/(1 + \eps^2)} 
Z_0^{2 \eps}} \geq 
\frac2{1 + \left(\rho(x/t)/ M_\eps
\right)^{-2 \eps/(1 + \eps^2)}}\\
& \geq \left( \rho(x/t) / M_\eps 
\right)^{\eps'}. 
\endaligned
$$
This establishes \eqref{52}.

Now, by expressing $W_l, Z_l, W_r, Z_r$ in terms of the initial density and velocity, one see that the constant $M_\eps$ can be taken to be independent of $\eps$. In fact, from the definition \eqref{modifiedRinvariants}, 
$$ W_0 = \rho_0^{1/(1 + \eps^2)} 
\left( \frac{1 + \eps v_0}{1 - \eps v_0}
\right)^{1/(2\eps)}, \, Z_0  = 
\rho_0^{1/(1 + \eps^2)} \left( 
\frac{1 + \eps v_0}{1 - \eps v_0}
\right)^{-1/(2\eps)}, $$ 
so 
$$\aligned 
W_0^{1 + \eps^2} & = 
\rho_0 \left( \frac{1 + \eps v_0}{1 
- \eps v_0}\right)^{(1 + \eps^2)/(2\eps)}\\
& = \rho_0 e^A \mbox{ with } A = \frac{1 + \eps^2}{2\eps} \log \left( \frac{1 + \eps v_0}{1 - \eps v_0}\right).\endaligned $$
By taking $\varphi(x)=\frac1{2x}\log \left( \frac{1+x}{1-x}\right)\, (x\not=0)$ with $\varphi (0) = 0$,
$$\aligned
A & = (1 + \eps^2) v_0 \, \varphi (\eps v_0)\\
& = (1 + \eps^2) v_0 \, \frac{\varphi''(0)}2 (\eps v_0)^2\left(1 + o(|\eps v_0|)\right) \endaligned $$
 and  $\varphi'(0) = 0,\, \varphi''(0) = 8/3$.
Hence
$$W_0^{1 + \eps^2} = \rho_0 \exp \left((1 + \eps^2) v_0 \, \frac{\varphi''(0)}2 (\eps v_0)^2\left(1 + o(|\eps v_0|)\right)\right)$$
Similarly, we can show 
$Z_0^{1 + \eps^2} = \rho_0 \exp \left((1 + \eps^2) v_0 \, \frac{\varphi''(0)}2 (\eps v_0)^2\left(1 + o(|\eps v_0|)\right)\right)$.
Therefore $M_\eps$ can be expressed only in terms of $(\rho_l,v_l)$ and 
$(\rho_r,v_r)$.
\end{proof}

Consider a family of cartesian discretizations of the spacetime $\RR_+ \times \RR$, based on
a time length $\tau$ and a space length $h$, where the ration $\tau/h$ is kept fixed while $h \to 0$.
Set $t_n := n \tau$ ($n$ being a positive integer) and $x_j := j \, h$ ($j$ being an integer).    
The Lax-Friedrichs scheme allows us to construct approximate solutions $\rho^h, v^h:\RR_+ \times \RR \to \RR$ 
to the relativistic Euler equations, which, for any pair $(n,j)$ with $n+j$ even,
is constant equal to $\rho_j^n, v_j^n$ in every slab $[t_n,t_{n+1}) \times I_j := [t_n,t_{n+1}) \times (x_{j-1}, x_{j+1})$.
The initial data $\rho_0, v_0$ is averaged over each initial cell
$$
(\rho_j^0, v_j^0) := \frac1{2h}\int_{I_j} (\rho_0, v_0)(x) \, dx, 
\quad j \text{ even.} 
$$
Given a piecewise constant approximation at a given time $t_n$, 
we solve a Riemann problem in the neighborhood of each point $x_{j+1}$, 
for $j$ such that $n+j$ is even, and we then average the solution at time $t=t_{n+1}-0$ 
over the intervals $ (x_{j}, x_{j+2})$: 
$$
(\rho_{j+1}^{n+1}, v_{j+1}^{n+1}) :=\frac1{2h} \int_{I_{j+1}} (\rho^h, v^h)(t_{n+1}-0,x) \, dx, \qquad n+j \text{ even}. 
$$

Observe that, if $\eps \in (0,1)$ (as well as for $\eps=1$), the characteristic speeds 
$\lambda_1, \lambda_2$ remain {\sl bounded globally,} 
even when the velocity approaches the light speed. Indeed, all wave speeds in the problem 
under consideration are bounded. 
(Note that a quite different situation is met with the non-relativistic model corresponding 
to $\eps=0$, for which 
the characteristic speeds are unbounded.)  
To avoid any interaction in two neighboring Riemann problems, it is 
necessary to restrict the ratio $\tau/h$ by the CFL stability condition 
\be
\label{CFL}
{\tau \over h} \max \Big( {|v^h| +1 \over 1 - \eps^2 |v^h|}\Big) < {1 \over 2},
\ee
where the maximum is taken over all $(t,x)$. 
For instance, a sufficient condition is
\be
\label{strongCFL}
{\tau \over h} {1  \over \eps} < {1 \over 2},
\ee 
which clearly becomes more restrictive as $\eps$ approaches $0$.

We observe that both the Riemann problem and the projection step satisfy uniform stability estimates. 
Finally, from Lemma~\ref{511} we deduce the key estimate of this section: 

\begin{lemma}[Uniform a~priori bounds]
\label{bounds4} Provided that the CFL condition \eqref{CFL} holds, 
the approximations $(\rho^h, v^h)$ satisfy the tame condition 
\be 
\aligned 
& 0 \leq \rho^h(t,x) \leq M, 
\\
& 1 - \eps |v^h(t,x)| \geq \big( \rho^h(t,x) /M\big)^{\eps'}, 
\endaligned 
\label{AS.bounds2}
\ee 
for some constant $M>0$ depending solely on the initial data. 
\end{lemma} 

In addition to the uniform amplitude bound on $\rho^h, v^h$ one can also derive 
entropy dissipation bounds which imply that the entropy inequalities associated with 
any weak entropy are satisfied up to an error term vanishing in the distribution sense as $h \to 0$. 
We omit the details and refer to DiPerna (1983) and Hsu, Lin, and Makino (2004) for further details 
on Lax-Friedrichs approximations.


\section{Reduction of the Young measure}
\label{Young}

We rely on the theory of compensated compactness for nonlinear hyperbolic systems for which we 
refer to Tartar (1979 \& 1983), Murat (1978 \& 1981), and DiPerna (1983).  
In view of the uniform bounds derived in Lemma~\ref{bounds4} we can associate to the sequence 
$\rho^h, v^h$ a Young measure $\nu=\nu_{t,x}$ supported 
in a tame region $\Tcal(M)$ for some uniform constant $M>0$.   
By definition, for almost every point, $\nu_{t,x}$ is a probability measure in the variable $\rho,v$
which has compact support in the Riemann invariant variables $W,Z$. The Young measure allows us 
to compute the weak limit of any composite function of the sequence $\rho^h, v^h$, that is: 
\be 
\label{Young2}
g(\rho^h, v^h) \to \la \nu, g \ra  
\ee 
in the sense of distributions, for every function $g$ that is continuous in the Riemann invariants $W,Z$. 

Relying on standard arguments one can check that the entropy dissipation measure   
associated with a weak entropy pair $(\Ucal,\Fcal)$, 
$$
\partial_t \Ucal (\rho^h,v^h) + \partial_x \Fcal (\rho^h,v^h)
$$
belongs to a compact set of the Sobolev space $H^{-1}_{\text{loc}}$. 
By the div-curl lemma, this property implies that $\nu$ satisfies Tartar's commutation relation
\be
\label{Tartar1}
\la \nu, \Ucal_1\Fcal_2 - \Ucal_2\Fcal_1 \ra 
= \la \nu, \Ucal_1 \ra \la \nu, \Fcal_2 \ra - \la \nu, \Ucal_2 \ra \la \nu, \Fcal_1 \ra
\ee
for any two weak entropy pairs $(\Ucal_1, \Fcal_1)$, $(\Ucal_2, \Fcal_2)$. Plugging the entropy-entropy flux pairs 
given by \eqref{BP.allentropies1} and \eqref{BP.allentropies3} in the
 relation \eqref{Tartar1} and dropping the test-function $\psi$, 
we obtain for almost every $(t,x)$ and for all $s,s' \in \RR$
\be
\label{Tartar2}
\la \nu, \chi(s) \sigma(s') - \chi(s) \sigma(s')  \ra 
= 
\la \nu, \chi(s) \ra \la \nu, \sigma(s') \ra - \la \nu, \chi(s') \ra \, \la \nu, \sigma(s)  \ra.  
\ee
Equivalently, this identity holds with $\sigma$ replaced by 
$\sigmas = \sigma - u \, \chi$.

Our main result in this section is as follows. We fix a point $(t,x)$ 
(where the above relation holds) and consider the Young measure at that point.

\begin{theorem}[Young measure reduction for the relativistic Euler equations]
\label{reduction}
A probability measure $\nu$ with compact support in the $(W,Z)$-plane and satisfying \eqref{Tartar2} for all 
$s,s'$ is either a Dirac mass 
$$
\nu = \delta_{\rho,v}
$$
or else has its support included in the vacuum line 
$$
\supp \nu \subset \big\{\rho = 0 \big\} = \big\{ WZ =0 \big\}. 
$$
\end{theorem}

In other words, there exists a function $(\rho,v)=(\rho,v)(t,x)$ defined for almost every $(t,x)$ such that, almost everywhere, 
either $\rho(t,x) >0$ and $\nu_{t,x} = \delta_{(\rho,v)(t,x)}$ or else 
$\rho(t,x) = 0$ and $\nu_{t,x}$ is supported in the vacuum line (the velocity begin then irrelevant). 
This immediately implies, for instance, 
that $\rho^h$ as well as $\rho^h v^h$ converge to their limits $\rho$ and $\rho \, v$, respectively.  
By a (standard) property of consistency of the Lax-Friedrichs scheme the limit must be a weak solution 
of the relativistic Euler equations satisfying entropy inequalities for all weak entropy pairs. 
This completes the proof of our main result stated earlier in Theorem~\ref{AS-maintheorem}.

\begin{proof} We closely follow the method used by LeFloch and Shelukhin (2005). 
Given some $s_1 \in \RR$, it will be convenient to write $\chi_1 = \chi(\rho,v,s_1)$. 
Given $s_1,s_2,s_3 \in\RR$, we consider the corresponding identity \eqref{Tartar2} for the three pairs
$$
(s_1,s_2), (s_2,s_3), (s_3,s_1).
$$
We then multiply these identities by
$\la \nu, \chi_3 \ra, \la \nu, \chi_1 \ra, \la \nu, \chi_2 \ra$,
respectively, 
where $\chi_j = \chi (\rho,v,s_j)\, (j=1,2,3)$, and 
we add up these three identities.
Then, due to the symmetry of the expressions, the sum of the right-hand sides vanishes identically:
$$
\aligned
&\la \nu, \chi_3 \ra \, 
\big(\la \nu, \chi_1 \ra \la \nu, \sigmas_2 \ra - \la \nu, \chi_2 \ra \, \la \nu, \sigmas_1  \ra\big)
+ 
\la \nu, \chi_1 \ra \big(\la \nu, \chi_2 \ra \la \nu, \sigmas_3 \ra - \la \nu, \chi_3 \ra \, \la \nu, \sigmas_2  \ra\big)
\\
&+\la \nu, \chi_2 \ra \, 
\big(\la \nu, \chi_3 \ra \la \nu, \sigmas_1 \ra - \la \nu, \chi_1 \ra \, \la \nu, \sigmas_3  \ra\big) = 0, 
\endaligned
$$
where $\sigmas_j=\sigmas(\rho,v,s_j)$, 
whereas the sum of the left-hand side is
$$
\la \nu, \chi_3 \ra \la \nu, \chi_1 \sigmas_2 - \chi_2 \sigmas_1 \ra + \la \nu, \chi_1 \ra \la \chi_2 \sigmas_3 - \chi_3 \sigmas_2 \ra + \la \nu, \chi_2 \ra \la \nu, \chi_3 \sigmas_1 - \chi_1 \sigmas_3 \ra = 0.
$$

Then, by differentiating once in $s_2$ and in $s_3$ and by setting
 $\partial_j := \partial /\partial s_j$, we obtain 
\be
\label{Tartar3}
\aligned
& \la \nu, \partial_3 \chi_3 \ra \la \nu, \chi_1 \partial_2 \sigmas_2 - \partial_2 \chi_2 \sigmas_1 \ra + \la \nu, \partial_2 \chi_2 \ra \la \nu, \partial_3 \chi_3 \sigmas_1 - \chi_1 \partial_3 \sigmas_3 \ra\\
= & - \la \nu, \chi_1 \ra \la \partial_2 \chi_2 \partial_3 \sigmas_3 - \partial_3 \chi_3 \partial_2 \sigmas_2 \ra,
\endaligned
\ee
which is valid in the sense of distributions in $\RR^3$. 
In view of \eqref{BP.singularities4a} and \eqref{BP.singularities4b}, 
by setting $X^{(i)}_j := X^{(i)}|_{s=s_j}$ we find, on one hand, 
\be
\label{commutationLH1}
\aligned
\chi_1\del_2 \sigmas_2 - \del_2 \chi_2 \sigmas_1 & = \big( \chi_1 X^{(3)} + \sigmas_1 X^{(1)} \big) \delta_{s_2 = w} + \big( \chi_1 X^{(3)} - \sigmas_1 X^{(1)} \big) \delta_{s_2 = z}\\
& \quad  + \big( \chi_1 X^{(4)}_2 - \sigmas_1 X^{(2)}_2 \big) \one_{E_2}
\endaligned
\ee
where $E_j := \big\{ (w-s_j)(z-s_j)<0 \big\}$ 
and
\be\label{commutationLH2}\aligned
\del_3 \chi_3 \sigmas_1 - \chi_1 \del_3 \sigmas_3 & = \big( - \sigmas_1 X^{(1)} - \chi_1 X^{(3)} \big) \delta_{s_3 = w} + \big( \sigmas_1 X^{(1)} - \chi_1 X^{(3)} \big) \delta_{s_3 = z}\\
& \quad  + \big( \sigmas_1 X^{(2)}_3 - \chi_1 X^{(4)}_3 \big) \one_{E_3}
\endaligned
\ee
and, on the other hand, 
\be
\label{commutationRH1}
\aligned
&\del_3 \chi_3 \del_2 \sigmas_2 - \del_2 \chi_2 \del_3 \sigmas_3\\
& = 2  X^{(1)} X^{(3)} \big( \delta_{s_2 = w} \delta_{s_3 = z} - \delta_{s_2 = z} \delta_{s_3 = w}  \big)\\
& \quad + \Big( \big( X^{(2)}_3 X^{(3)} + X^{(1)} X^{(4)}_3 \big) \, \delta_{s_2 = w} 
  +  \big( X^{(2)}_3 X^{(3)} - X^{(1)} X^{(4)}_3 \big) \, \delta_{s_2 = z}
  \Big) \one_{E_3}\\
& \quad + \Big(  \big( - X^{(1)} X^{(4)}_2 - X^{(2)}_2 X^{(3)} \big) \, \delta_{s_3 = w} 
   + \big( X^{(1)} X^{(4)}_2 - X^{(2)}_2 X^{(3)} \big) \, \delta_{s_3 = z}  \Big) \,
    \one_{E_2}
\\
&  \quad + \big( X^{(2)}_3 X^{(4)}_2 - X^{(2)}_2 X^{(4)}_3 \big) \one_{E_2} \one_{E_3}. 
\endaligned
\ee

In view of the formulas \eqref{commutationLH1} and \eqref{commutationLH2}, the right-hand side of \eqref{Tartar3} contains products of functions with bounded variation and Dirac masses, plus regular terms. 
Such products were already discussed by Dal~Maso, LeFloch, and Murat (1995). 
On the other hand, the right-hand side of \eqref{Tartar3} is more singular 
and involves also products of measures. 
Our calculations will show that the left-hand side of \eqref{Tartar3} tends to zero in the sense of 
distributions when $s_2, s_3 \to s_1$, while the right-hand side tends to a non-trivial limit.

We test the identity \eqref{Tartar3} with the following function of the variables $s_1, s_2, s_3$
\be\label{test}
\psi (s_1) \, \varphi^\kappa_2 (s_1 - s_2) \, \varphi^\kappa_3 (s_1 - s_3) 
:= \psi (s_1) \, \frac1{\kappa^2} \varphi_2 \big( (s_1 - s_2)/\kappa \big) 
    \, \varphi_3 \big( (s_1 - s_3)/\kappa \big), 
\ee 
where $\kappa$ is a small parameter, 
$\psi$ is a smooth and compactly supported function, 
and $\varphi_j : \RR \rightarrow \RR \, (j=2,3)$ are mollifiers such that
$$
\varphi_j  \geq 0, \quad 
\int_\RR \varphi_j (s_j) \, ds_j = 1, 
\quad 
\supp \varphi_j \subseteq (-1,1).
$$

We consider first the right-hand side of \eqref{Tartar3}. 
Applying the test-function \eqref{test} to the distribution $\la \nu, \chi_1 \ra \la \partial_2 \chi_2 \partial_3 \sigmas_3 - \partial_3 \chi_3 \partial_2 \sigmas_2 \ra$, 
we obtain the integral term 
\be
\label{commutationRH2}
\int_{\RR^3}\la \nu, \chi_1 \ra \la \partial_2 \chi_2 \partial_3 \sigmas_3 - \partial_3 \chi_3 \partial_2 \sigmas_2 \ra
\psi (s_1) \varphi^\kappa_2 (s_1 - s_2) \varphi^\kappa_3 (s_1 - s_3)\, ds_1 ds_2 ds_3, 
\ee
which we decompose as a sum $\sum\limits_{i=1}^4 I_i^\kappa$.
Relying here on Lemma \ref{BP.singularities4}, we distinguish between products of Dirac measures
$$
I_1^\kappa := 2 \int_\RR \psi (s_1) \la \nu, \chi_1 \ra \left\la X^{(1)} X^{(3)} \left(  \varphi_2^\kappa (s_1 - w) \varphi_3^\kappa (s_1 - z) - \varphi_2^\kappa (s_1 - z) \varphi_3^\kappa (s_1 - w)\right) \right\ra \, ds_1,
$$
products of Dirac measures by functions with bounded variation
$$\aligned
I_2^\kappa := & \int_\RR \psi (s_1) \la \nu, \chi_1 \ra \big\la  \varphi_2^\kappa (s_1 - w) \int_\RR \left( X^{(2)}_3 X^{(3)} + X^{(1)} X^{(4)}_3 \right) \\
& \hskip5.cm \varphi_3^\kappa (s_1 - s_3) \one_{E_3} \, ds_3 \big\ra \, ds_1\\
& - \int_\RR \psi (s_1) \la \nu, \chi_1 \ra\big\la \varphi_3^\kappa (s_1 - w) \int_\RR \left( X^{(2)}_2 X^{(3)} +  X^{(1)} X^{(4)}_2\right)
\\
& \hskip5.cm \,\varphi_2^\kappa (s_1 - s_2) \one_{E_2} \, ds_2 \big\ra \, ds_1\\
=: & I_{2,1}^\kappa - I_{2,2}^\kappa, 
\endaligned 
$$
$$
\aligned 
I_3^\kappa := & \int_\RR \psi (s_1) \la \nu, \chi_1 \ra \big\la \varphi_2^\kappa (s_1 - z) \int_\RR \left( X^{(2)}_3 X^{(3)} - X^{(1)} X^{(4)}_3 \right)\\
& \hskip5.cm \,\varphi_3^\kappa (s_1 - s_3) \one_{E_3} \, ds_3 \big\ra \, ds_1\\
& - \int_\RR \psi (s_1) \la \nu, \chi_1 \ra \big\la \varphi_3^\kappa (s_1 - z) \int_\RR \left(  X^{(2)}_2 X^{(3)} - X^{(1)} X^{(4)}_2 \right)\\
& \hskip5.cm \, \varphi_2^\kappa (s_1 - s_2) \one_{E_2} \, ds_2 \big\ra \, ds_1. 
\endaligned
$$
The remainder
$$
\aligned
I_4^\kappa := & \int_{\RR^3} \psi (s_1) \la \nu, \chi_1 \ra \big\la \left( X^{(2)}_3 X^{(4)}_2 - X^{(2)}_2 X^{(4)}_3 \right) \one_{E_2} \one_{E_3}\big\ra\\
& \hskip5.cm \varphi_2^\kappa (s_1 - s_2) \varphi_3^\kappa (s_1 - s_3) \, ds_1 ds_2 ds_3
\endaligned
$$
involve functions only. 

By a change of variable we see that the integral term 
$$
\aligned
I_1^\kappa = & \frac2\kappa \int_{W,Z}\int_{W',Z'}\int_\RR \psi( \kappa y + w) \chi(\rho,v,\kappa y + w) X^{(1)} X^{(3)}\\
& \left( \varphi_2^\kappa (y) \varphi_3^\kappa (y + (w - z)/\kappa) - \varphi_2^\kappa (y + (w - z)/\kappa) \varphi_3^\kappa (y) \right) \, dy d\nu d\nu'
\endaligned
$$
tends to zero, i.e. $I_1^\kappa \to 0$. The same is true for the smoothest term $I_4^\kappa$, in fact
$$\aligned
I_4^\kappa = & \int_{W,Z}\int_{W',Z'}\int_\RR \Big( \int_z^w X^{(2)} (\rho,v,s_3) \varphi_3^\kappa (s_1 - s_3) \, ds_3 \int_z^w X^{(4)} (\rho,v,s_2)\\
& \, \varphi_2^\kappa (s_1 - s_2) \, ds_2 - \int_z^w X^{(2)}(\rho,v,s_2) \varphi_2^\kappa (s_1 - s_2) \, ds_2 \int_z^w X^{(4)} (\rho,v,s_3)\\
& \, \varphi_3^\kappa (s_1 - s_3) \, ds_3 \Big)  \psi (s_1) \chi (\rho,v,s_1) \, ds_1 d\nu d\nu',
\endaligned
$$
which clearly tends to
$$\aligned
&\int_{W,Z}\int_{W',Z'}\int_\RR \left( X^{(2)}(\rho,v,s_1) X^{(4)} (\rho,v,s_1) - X^{(2)}(\rho,v,s_1) X^{(4)}(\rho,v,s_1) \right)\\
& \hskip7.cm \, \psi (s_1) \chi (\rho,v,s_1) \, ds_1 d\nu d\nu' =0.
\endaligned
$$

Next, let us consider the term $I_2^\kappa = I_{2,1}^\kappa - I_{2,2}^\kappa$ in \eqref{commutationRH2}. 
By defining
$$
\aligned 
& Q^- : =X^{(2)} X^{(3)} + X^{(1)} X^{(4)}, 
\\
& Q^-_j := Q^- (\rho,v,s_j), \qquad 
\chi' := \chi(\rho',v',s), \chi_j' := \chi (\rho',v's_j), 
\endaligned
$$
we can write 
$$
I_{2,1}^\kappa = -\int_{W,Z} \int_{W',Z'} \int_\RR \psi(s_1) \chi_1'  \varphi_2^\kappa (s_1 - w) \int_w^z Q^-_3 \varphi_3^\kappa (s_1 - s_3) \, ds_3 ds_1 d\nu d\nu'. 
$$

\begin{lemma}
\label{product1}
Let $f, F : \RR \rightarrow \RR$ be continuous functions. Then, for 
every interval $[a,b], [a',b'] \subseteq \RR$, the integral term 
$$
I^\kappa (a,b,a',b') := \int_{a'}^{b'}  f(s_1) \varphi_2^\kappa (s_1 - a) \int_a^b F(s_3) \varphi_3^\kappa (s_1 - s_3) \, ds_3 ds_1
$$
converges (when $\kappa \to 0$) toward 
$$
f(a) \, F(a) \, \left( A_{2,3} \, \one_{a' < a <b'} + B_{2,3} \, \one_{a=a'} + C_{2,3} \, \one_{a=b'} \right),
$$
where the coefficients $A_{2,3} := B_{2,3} + C_{2,3}$ and $B_{2,3}$ and $C_{2,3}$ 
depend on the mollifying functions as follows: 
$$
\aligned
B_{2,3} : = & \int_0^\infty \int_{-\infty}^{y_1} 
\varphi_2 (y_1) \varphi_3 (y_3) \, dy_3 dy_1,
\\
C_{2,3} := & \int_{-\infty}^0 \int_{- \infty}^{y_1} 
\varphi_2 (y_1) \varphi_3 (y_3) \, dy_3 dy_1.
\endaligned
$$
\end{lemma}

Formally the integral $I^\kappa$ has the limit 
$$
I (a,b,a',b') := \int_{a'}^{b'} f(s_1) \, \delta_{s_1 = a} \, \int_a^b F(s_3) \delta_{s_3 = s_1} \, ds_3 ds_1.
$$

\begin{proof}
Making first the change of variables $s_3 = s_1 - \kappa y_3$
and then $s_1 = \kappa y_1 + a$, we can write
$$
I^\kappa  = -\int_{(a'-a)/\kappa}^{(b'-a)/\kappa} f(\kappa y_1 + a) \varphi_2 (y_1)
 \int_{y_1}^{y_1-(b-a)/\kappa} F\left(\kappa(y_1 -  y_3)+a\right) \varphi_3 (y_3) 
\, dy_3 dy_1.
$$
Clearly, we have $I^\kappa \to 0$ when $a<a'$ or $a>b'$.

Now, if $a=a'$, we can write
$$\aligned 
I^\kappa & = - \int_0^{(b'-a)/\kappa} f(\kappa y_1
+ a) \varphi_2 (y_1) \int_{y_1}^{y_1 -
(b-a)/\kappa} F \left( \kappa(y_1 - y_3)+a) \right) 
\varphi_3(y_3) \, dy_3 dy_1\\
& \rightarrow f(a) F(a) \int_0^\infty  
\varphi_2(y_1) \int_{-\infty}^{y_1}
\varphi_3 (y_3) \, dy_3dy_1
\endaligned
$$
as $\kappa \to 0$. The other values of $a$ can
be studied by the same argument and this
completes the proof of Lemma \ref{product1}.
\end{proof}

Applying Lemma \ref{product1}, we find that $I_{2,1}^\kappa$ tends to
$$
\aligned
& -\int_{W,Z}\int_{W',Z'} \psi (w)\chi' (w) Q^- (w) \big(  A_{2,3} \one_{-\infty < w<\infty} + B_{2,3} \one_{w = - \infty} + C_{2,3} \one_{w=\infty} \big)\, d\nu d\nu'\\
& = - A_{2,3} \int_{W,Z}\int_{W',Z'} \psi (w)\chi' (w) Q^- (w) \, d\nu d\nu'
\endaligned
$$
and that $I_{2,2}^\kappa$ tends to
$$ 
- A_{3,2} \int_{W,Z}\int_{W',Z'} \psi (w)\chi' (w) Q^- (w) \, d\nu d\nu'. 
$$

We conclude that the limit of $I_2^\kappa$ is 
$$
A \int_{W,Z}\int_{W',Z'} \psi (w) \chi'(w) Q^- (w) \, d\nu d\nu',
$$
where
$$
\aligned
A := & A_{3,2} - A_{2,3}\\
= & \int_\RR \int_{-\infty}^{y_1} \left(
\varphi_2 (y_3) \varphi_3 (y_1)
 - \varphi_2 (y_1) \varphi_3 (y_3)\right)
\, dy_3 dy_1.
\endaligned
$$

We can determine similarly that the limit of $I_3^\kappa$ which is found to be 
$$ 
- A \int_{W,Z}\int_{W'.Z'} \psi (z) \chi'(z) Q^+(z) \, d\nu d\nu',
$$
where
$$
Q^+ := X^{(2)} X^{(3)} - X^{(1)} X^{(4)}.
$$
In conclusion, we have identified the limit of the term \eqref{commutationRH2}, as 
$$
A \left(\int_{W,Z}\int_{W',Z'} \psi (w) \chi'(w) Q^- (w) \, d\nu d\nu' -  \int_{W,Z}\int_{W'.Z'} \psi (z) \chi'(z) Q^+(z) \, d\nu d\nu' \right).
$$

One can also check (LeFloch and Shelukhin, 2005) that the distribution 
$$
\la \nu, \partial_3 \chi_3 \ra \la \nu, \chi_1 \partial_2 \sigmas_2 - \partial_2 \chi_2 \sigmas_1 \ra + \la \nu, \partial_2 \chi_2 \ra \la \nu, \partial_3 \chi_3 \sigmas_1 - \chi_1 \partial_3 \sigmas_3 \ra
$$
(which is the left-hand side of \eqref{Tartar3} applied by the test-function \eqref{test} 
tends to zero as $\kappa \to 0$.

Since the molifying functions $\varphi_2$ and $\varphi_3$ can always be chosen such that $A \not=0$, 
we have reached the following conclusion 
\be
\label{conclu} 
\aligned
& \int_{W,Z}\int_{W',Z'} \psi (w) \chi'(w) Q^- (w) \, d\nu d\nu' 
-  \int_{W,Z}\int_{W'.Z'} \psi (z) \chi'(z) Q^+(z) \, d\nu d\nu'
\\ 
& = \left\la \nu \otimes \nu', \psi (w) \chi'(w) Q^- (w) 
-   \psi (z) \chi'(z) Q^+(z)  \right\ra=0. 
\endaligned
\ee
By Lemma \ref{BP.singularities4} we have
$$\aligned
Q^-(w) = & X^{(2)} (0, w - z) X^{(3)} (w-z)
 +X^{(1)} (w - z) X^{(4)} (0, w - z),
 \\
 Q^+ (z) = & X^{(2)} (-w +z,0) X^{(3)}(w-z) 
 - X^{(1)}(w-z) X^{(4)} (-w+z,0). 
\endaligned
$$

Observe that the test-function $\psi$ is arbitrary and 
can be used to localize the equation \eqref{conclu}. In turn, 
we end up with a necessary condition satisfied by the Young measure $\nu$:  
$$
\la \nu \otimes \nu', \Xi \, g \ra = 0.  
$$
Here, 
$\nu \otimes \nu'$ denotes the tensor product of $\nu$ and $\nu'$ (another copy of the Young measure $\nu$), 
while $g$ is some (non-negative) characteristic function, 
while the function $\Xi$ is defined by
$$ 
\Xi(\rho) := \lim_{z \to 0, wz <0}\chi(w,z) \, \left( \sigmas_w(w,z) + \sigmas_z(w,z)\right) 
             + \left(\chi_w(w,z) + \chi_z(w,z)\right) \, \sigmas(w,z), 
$$
$\rho$ and $w$ being related by $w=2 R = 2 (\ln \rho)/(1+\eps^2)$. 
An analogous identity holds with $\lim\limits_{z \to 0, wz <0}$ replaced by $\lim\limits_{w \to 0, wz <0}$. 
For convenience, we write simply $\chi (w,0)$ instead of $\lim\limits_{z \to 0, wz <0} \chi (w,z)$, etc. 
This allows to rewrite the expression of $\Xi$ as 
\be
\label{Xi}
\aligned
\Xi(w) & = - \left(\chi (w,0,0)  \sigmas_s (w,0,0)   + \chi_s (w,0,0) \sigmas (w,0,0) \right)\\
& = \chi (w,0)  \left(\sigmas_w (w,0)   +  \sigmas_z (w,0)\right) + \left( \chi_w (w,0) + \chi_z (w,0) \right) \sigmas (w,0). 
\endaligned
\ee

{\sl Provided} the coefficient $\Xi$ keeps a constant sign, the above condition implies 
(LeFloch and Shelukhin, 2005) 
that the support of $\nu$ either is a single point, or else has its support concentrated 
where $WZ =0$,
which is nothing but the vacuum line. 
Recall that the relativistic equations are automatically satisfied if the density vanishes identically.

It remains to establish that $\Xi(\rho) < 0$. 
More precisely, we only need this to hold for sufficiently small $\rho$,
since by using the scaling invariance property of the relativistic Euler equations
(Lemma \ref{invariance}) we can always ensure that the range of the Lax-Friedrichs approximations 
and therefore the support of the Young measure is included in a neighborhood of the vacuum. 
We will actually prove that 
\be
\label{inequal1}
\Xi(\rho) \leq - {1 \over 2} \rho^{2\alpha} < 0, \qquad 
0< \rho  \ll 1. 
\ee

We set $\Xi(w) =: \lim\limits_{z \to 0, wz <0}\chi^2(w,z) \, \Omega(w)$, and we observe that $\lim\limits_{z \to 0, wz <0}\chi(w,z)$ behaves like $\rho^\alpha$ 
near the vacuum and, therefore, we need to check that 
\be
\label{ineql3}
\aligned 
\Omega(w) 
& = \left( - \frac12 + \frac{w}8 \right) + O (\eps^2)\\
& < 0. 
\endaligned 
\ee 
The term $O(\eps^2)$ should be bounded by a constant times $\eps^2$, 
uniformly for all large (negative) $w$. This will show that
$\Omega$ has a fixed sign for all sufficiently negative values of $w$, 
and this will indeed complete the proof of the theorem.

\

It remains to determine the sign of the function $\Xi$. 
For the sake of comparison, we recall first the relevant formula when $\eps=0$: 
$$
\lambda_1^0 (w,0) = \frac{w}2 - 1, \qquad \lambda^0_2 (w,0) = \frac{w}2 + 1,
$$
and
$$
\aligned 
& \chi^0(w,0) = e^{w/4}, 
\\
& \chi_w^0 (w,0) = \frac14 e^{w/4}, \, \chi_z (w,0) = - \frac14 ( 1 + w/4) e^{w/4},
\\
&  (\chi_w^0 + \chi_z^0)(w,0) = -  \frac{w}{16} \, e^{w/4},
\endaligned 
$$ 
and 
$$
\aligned 
& \sigma^0(w,0) = \left(\frac{w}2 - 1\right) e^{w/4}, \, \sigmas (w,0) = - e^{w/4}, 
\\
& (\sigma_w^0 + \sigma_z^0)(w,0) = \left(\frac12 + \frac{w}{16} - \frac{w^2}{32}\right),
\\
& (\sigma^{\sharp,0}_w + \sigma^{\sharp,0}_z )(w,0) = \left( - \frac12 + \frac{w}{16} \right) e^{w/4}. 
\endaligned 
$$ 
Hence, in the {\sl non-relativistic case} we do have 
$$
\Xi^0(w) = \left( - \frac12 + \frac18 w \right) e^{w/2}.  
$$

For the general case $\eps \neq 0$ we have already 
determined the traces of the entropy kernel along the boundary $z=0$ while keeping $wz<0$
(see \eqref{chi.value0}, \eqref{chi.value1a} and \eqref{BP.singularities3b}): 
\be
\label{chi1}
\aligned 
\chi(w,0) & = \frac{1 + \eps + (1 - \eps) e^{\eps w}}2 \, e^{ (1 - \eps)^2 w / 4},
\\
\chi_w(w,0) & = {1 - \eps^2\over 8} \big(1 - \eps + (1 + \eps) e^{\eps w} \big) \, e^{ (1 - \eps)^2 w / 4}
\\
& =  {1 - \eps^2\over 4} {1 - \eps + (1 + \eps) e^{\eps w}  \over 1 + \eps + (1 - \eps) e^{\eps w}  } \,  \chi(w,0)= \abar (w)\chi(w,0),
\\
\chi_z(w,0)  
& = - {1 - \eps^2 \over 4} \left( 1 + (1-\eps^2) w/4 \right) \chi(w,0). 
\endaligned 
\ee

On the other hand, for the entropy kernel, 
from \eqref{BP.entropypair2}, \eqref{speeds} and \eqref{chi.value1a} it follows that 
$$
\aligned 
\sigma_w(w,0) 
& = \lambda_2(w,0) \, \chi_w(w,0), 
\\
& =  -{1 \over \eps} {1 - \eps - (1+\eps) e^{\eps w} \over 1- \eps + (1+\eps) e^{\eps w}}
     \, {1 - \eps^2\over 4} \frac{1 - \eps + (1 + \eps) e^{\eps w}}2 \, e^{ (1 - \eps)^2 w / 4}\\
& = -  (1 - \eps^2) \,  \frac{1 - \eps - (1+\eps) e^{\eps w}}{8 \eps} \, e^{ (1 - \eps)^2 w / 4}, 
\endaligned 
$$
that is
$$
\aligned  
\sigma(w,0) 
& = \int_{-\infty}^w \sigma_w(w',0) dw' 
\\
& =  - \int_{-\infty}^w (1 - \eps^2)   \frac{1 - \eps - (1+\eps) e^{\eps w'}}{8 \eps} \, e^{ (1 - \eps)^2 w' / 4} \, dw'. 
\endaligned 
$$
So, we find 
$$
\aligned
\sigma(w,0) & = -\frac{1 + \eps - (1 - \eps) e^{\eps w}}{2\eps} e^{(1 - \eps)^2 w/4}\\
& = -{1 \over \eps} {1+\eps - (1-\eps) e^{\eps w} \over 1+ \eps + (1-\eps) e^{\eps w}} \, \chi(w,0) = \lambda_1 (w,0) \chi (w,0)\\ 
\endaligned
$$
and thus in view of \eqref{lorentz3} 
\be
\label{sig1}
\aligned 
\sigmas(w,0) & = \left( \lambda_1(w,0) - \frac{w}2 \right) \chi(w,0), 
\\
\sigmas_w (w,0) & = \left( \lambda_2 (w,0) - \frac{w}2 \right) \chi_w (w,0) - \frac12 \chi (w,0).  
\endaligned 
\ee
Furthermore, from \eqref{lorentz3} we get 
\be
\label{sig2}
\sigmas_z(w,0)  = \left(\lambda_1(w,0)-\frac{w}2\right) \chi_z (w,0) - \frac12 \chi(w,0). 
\ee

Finally, we are in a position to compute the quantity of interest.  
In view of \eqref{chi1}, \eqref{sig1}, and \eqref{sig2}, we can obtain 
$$ 
\aligned 
& \Xi^{(w)}(\rho) := \chi(w,0) \, \sigmas_w(w,0)  + \chi_w(w,0) \, \sigmas(w,0), 
\\
& \Xi^{(z)}(\rho) := \chi(w,0) \, \sigmas_z(w,0) + \chi_z(w,0) \, \sigmas(w,0), 
\endaligned 
$$
as follows: 
$$
\aligned
\Xi^{(z)}(\rho)
& = \chi (w,0) \, \left( \left( \lambda_1(w,0) - \frac{w}2 \right) \chi_z (w,0) - \frac12 \chi(w,0) \right)
\\
&\hspace{0.2cm} + \chi_z (w,0) \left( \lambda_1(w,0)  - \frac{w}2\right) \chi(w,0) 
\\
& = \chi (w,0) \left( 2 \left( \lambda_1 (w,0) - \frac{w}2 \right) \chi_z (w,0) - \frac12 \chi (w,0) \right)
\\
& =  - \frac12 \chi^2(w,0) \left( (1 - \eps^2)\left( 1 + (1-\eps^2)w/4 \right) \left( \lambda_1(w,0) - \frac{w}2 \right) +1 \right)
\endaligned
$$
and 
$$
\aligned
\Xi^{(w)}(\rho)
& = \chi(w,0) \left(\left( \lambda_2 (w,0) - \frac{w}2\right) \chi_w (w,0) - \frac12 \chi(w,0)\right)
\\
& \quad + \chi_w(w,0) \, \left( \lambda_1(w,0) - \frac{w}2 \right) \chi(w,0) 
\\
& = \chi(w,0) \left( \left(  \lambda_1 (w,0) + \lambda_2 (w,0) - w \right) \chi_w (w,0) - \frac12 \chi (w,0)\right)
\\
& = \chi(w,0)^2 \left(\left(  \lambda_1(w,0) + \lambda_2 (w,0) - w \right) \abar (w) - \frac12 \right). 
\endaligned 
$$ 

Thus, we conclude that $\Xi (w)= \chi^2(w,0) \, \Omega(w)$, where 
\be
\label{Key}
\aligned 
\Omega(w) := 
&  -\frac12(1 - \eps^2)\left( 1 + (1-\eps^2)w/4 \right) \left( \lambda_1(w,0) - \frac{w}2 \right)\\
& + \big( \lambda_1(w,0) + \lambda_2 (w,0) - w \big) \abar (w) - 1.
\endaligned 
\ee
We easily see that
$$
\aligned 
& \lambda_1 (w,0) = w/2 -  1 + O (\eps^2), 
\\
& \lambda_2 (w,0) = w/2 + 1 + O (\eps^2), 
\\
& \abar (w) = - 1/4 + O(\eps^2), 
\endaligned 
$$
thus \eqref{ineql3} holds and $\Xi$ vanishes only at the vacuum $\rho = 0$.
\end{proof} 

\

The formulas derived in the present paper converges formally to the ones of the non-relativistic 
case, as we now show by using the notation introduced by LeFloch and Shelukhin (2005).  
Just before the identity $(5.7)$ on p.~424 of that paper,  
the function $D(R) := Q^- (w)$ is defined while 
$Q^- (w):= G^\chi - G^h$ is introduced on p.~420.
The terms $G^\chi$ and $G^h$ are introduced in Theorem~4.6 on p.~414, 
that is: 
$$
\aligned
G^\chi (R,u) & = - 2 |R|  f' (0) e^{R/2}, \qquad u \leq - |R|, 
\\
G^h (R,u) & = e^{R/2} \big( 2R + \frac12\big), \qquad |u| \geq |R|.
\endaligned
$$
The calculation in the proof of Theorem~4.6 should be modified, as follows. 

First of all, given a test function $\varphi = \varphi (s)$, we can write
$$\aligned
\langle \chi', \varphi \rangle & = - \int_\RR \chi (s) \varphi'(s) \, ds
 = - e^{R/2} \int_{u - |R|}^{u + |R|} \varphi'(s) f(| u - s|^2 - R^2) \, ds\\
& = - e^{R/2} \left( \varphi ( u + |R|) - \varphi (u - |R|) \right)\\
& \quad - e^{R/2} \int_{u - |R|}^{u + |R|} 2 \varphi(s) f'(|u - s|^2 - R^2) (u-s)\, ds,
\endaligned
$$
which yields
$$\chi_s = e^{R/2} \left( \delta_{s= u - |R|} - \delta_{s = u + |R|} \right) + G^\chi (R, u - s) \one_{|u - s| < |R|},
$$
where, for all $|v| \leq |R|$,
$$
G^\chi (R,v) = - 2 e^{R/2} v f'(v^2 - R^2).
$$
Hence, we find 
$$
G^\chi (R,u) = 2 |R| f'(0) e^{R/2}, \qquad u \leq - |R|, 
$$
and this formula contains a plus sign, instead of a minus sign as stated originally.

Second, for the expression of the derivative of the entropy flux $h_s$,  
a term $f'(0)$ should be added, as follows. We write 
$$
\aligned
h_s
& = \sgn (u - s) e^{- |u - s|/2} \Big( \delta_{s = u - |R|} - \delta_{s = u + |R|} + \frac12 \sgn (u - s) \one_{|u - s| < |R|}\Big)\\
&  \quad - 2 \int_{- (|R| \vee |u - s|)}^{- |u - s|} \frac\del{\del s}\left( ( u - s) e^{r/2} f'(|u - s|^2 - r^2) \right) \, dr\\ 
&  \quad 
+ 2 f'(0) e^{- |u - s|/2} (u - s) \sgn (u - s)  - 2 e^{- |u - s|/2} \one_{|u - s| \geq |R|} f'(0) |u - s|
\\
& = e^{R/2} \left( \delta_{s= u - |R|} + \delta_{s= u + |R|}\right) + \frac12 e^{- |u - s|/2} \one_{|u - s| < |R|}\\
& \quad + 2 \int_{- ( |R| \vee |u - s|)}^{- |u - s|} \left( e^{r/2} f'( |u - s|^2 - r^2) + 2 e^{r/2} |u - s|^2 f''(|u - s|^2 - r^2) \right) \, dr\\
& \quad + 2 f'(0) e^{- |u - s|/2} |u - s| \one_{|u - s| < |R|}
\endaligned
$$
and, therefore,
$$
h_s = e^{R/2} \left( \delta_{s= u - |R|} + \delta_{s= u + |R|}\right) + G^h (R, u-s) \one_{|u - s| <|R|}
$$
with
$$\aligned
G^h (R,v) : & = e^{-|v|/2} \, \big( \frac12 + 2 f'(0) |v| \big)
\\
& = e^{-|v|/2} \, \big( \frac12 - \frac18 |v| \big).
\endaligned
$$
Hence, we find 
$$
G^h(R,u) = e^{R/2} \, \big( \frac12 + \frac{R}8 \big), \qquad u\leq - |R|.
$$

In conclusion, for small $\rho$ we have $R <0$ and 
$$\aligned
Q^- (w) &= - 2 R \, \big(- \frac1{16} \big) \, e^{R/2} - e^{R/2} \, \big(\frac12 + \frac{R}8\big)
\\
& = e^{R/2} \, \big(- \frac12 + \frac{R}4\big) = e^{w/4}\, \big(- \frac12 + \frac{w}8\big), 
\endaligned
$$
which is precisely the same expression as the limit $\eps \to 0$ of the expression 
\eqref{ineql3} obtained in the present paper. 


\section*{Acknowledgments} 

The first author (PLF) is very thankful to the organizers (P.T. Chrusciel, H. Friedrichs, P. Tod) 
of the Semester Program
``Global Problems in Mathematical Relativity'' which took place at the Isaac Newton Institute of Mathematical Sciences 
(Cambridge, UK) and where this research was initiated. 

PLF was partially supported by the A.N.R. (Agence Nationale de la Recherche)
through the grant 06-2-134423
entitled {\sl ``Mathematical Methods in General Relativity''} (MATH-GR), and by the Centre National de la Recherche Scientifique (CNRS). The second author (MY) was supported by a Grant-in-Aid for Scientific Research
from the Japan Society for the Promotion of Science (JSPS). 


\end{document}